\documentclass{mcom-l}

\usepackage{rotating}
\usepackage{subfig}
\usepackage{amsmath}
\usepackage[english]{babel}
\usepackage{setspace}
\usepackage{amsfonts}
\usepackage{graphicx}
\usepackage{epstopdf}
\usepackage{algorithm}
\usepackage{algorithmic}
\usepackage{mathtools}
\usepackage{multirow}
\usepackage{booktabs} 
\usepackage{makecell}
\usepackage{url}
\usepackage{cleveref}

\Crefname{equation}{}{}

\newtheorem{theorem}{Theorem}[section]
\newtheorem{lemma}{Lemma}[section]
\theoremstyle{definition}

\theoremstyle{remark}

\def\mB{\mathcal{B}}
\def\mC{\mathcal{C}}
\def\mD{\mathcal{D}}

\def\mI{\mathcal{I}}

\def\mM{\mathcal{M}}

\def\mO{\mathcal{O}}

\def\mR{\mathcal{R}}
\def\mS{\mathcal{S}}

\def\md{\mathrm{d}}

\begin{document}

\title{Spectral Approximation to Fractional Integral Operators}

\author{Xiaolin Liu}
\address{School of Mathematical Sciences, University of Science and Technology of China, 96 Jinzhai Road, Hefei 230026, Anhui, China}
\email{xiaolin@mail.ustc.edu.cn}

\author{Kuan Xu}
\address{School of Mathematical Sciences, University of Science and Technology of China, 96 Jinzhai Road, Hefei 230026, Anhui, China}
\email{kuanxu@ustc.edu.cn}

\subjclass[2020]{34A08, 45E10, 26A33, 65R10, 47A58, 34K37, 44A35}

\keywords{fractional integral, fractional differentiation, operator approximation, spectral method, Chebyshev polynomials}

\begin{abstract}
We propose a fast and stable method for constructing matrix approximations to fractional integral operators applied to series in Chebyshev fractional polynomials. Based on a recurrence relation satisfied by the definite integrals of mapped Chebyshev polynomials with a fractional weight, the proposed method significantly outperforms existing approaches. Through numerical examples, we highlight the broad applicability of these matrix approximations, including the solution of boundary value problems for fractional integral and differential equations. Additional applications include fractional differential equation initial value problems and fractional eigenvalue problems.
\end{abstract}

\maketitle

\section{Introduction}\label{sec:intro}
This paper focuses on the fast and stable construction of spectral approximations to the fractional integral operator (FIO)
\begin{align*}
\mathcal{I}^{\mu}u(x) = \frac{1}{\Gamma(\mu)} \int_{-1}^{x}\frac{u(t)}{\left(x-t\right)^{1-\mu}} \md t,
\end{align*}
for $x \in [-1,1]$ and $\mu > 0$, where $\mu$ may be rational or irrational and $\Gamma(\cdot)$ denotes the gamma function. This operator can be deemed as a convolution of Volterra type with a weakly singular kernel and is a fundamental building block in both fractional integral equations (FIEs) and fractional differential equations (FDEs).

We work with a Chebyshev-based version of Jacobi fractional polynomials\footnote{The Jacobi fractional polynomials are also referred to as M\"{u}ntz Jacobi polynomials \cite{xu2}.} (JFPs) \cite{kaz,bhr,pu} defined as
\begin{align}
Q_n^{\alpha,\beta}(x) = \left(\frac{1+x}{2}\right)^{\alpha} T_n\left(2\left(\frac{1+x}{2}\right)^{\beta}-1\right), \label{jfp}
\end{align}
where $T_n$ is the $n$th Chebyshev polynomial of the first kind. It has been shown that the set $\{Q_n^{\alpha, \beta}(x)\}_{n=0}^{\infty}$ forms an orthogonal basis with respect to the weight function
\begin{align*}
w(x) = (1+x)^{\frac{\beta}{2}-1 - 2\alpha} \left(1 - 2^{-\beta}(1+x)^{\beta}\right)^{-\frac{1}{2}};
\end{align*}
see also \cite[\S 3.2]{bhr}. The parameter $\beta > 0$ is selected such that $\mu = k\beta$ for some $k \in \mathbb{N}^{+}$, while $\alpha > -1$ is typically chosen to reflect the singular behavior of the problem, for instance, to match the singularities in the right-hand side $f(x)$ of an FIE (see \Cref{fie} below). Notably, neither $\alpha$ nor $\beta$ is restricted to rational values.

It can be shown that $\{Q_n^{\alpha,\beta}(x)\}_{n=0}^{\infty}$ can be generated via Gram--Schmidt orthogonalization of $\{(\mI^{\mu})^j Q_0^{\alpha,\beta}(x)\}_{j=0}^{\infty}$, which span the Krylov space of $\mathcal{I}^{\mu}$ \cite[(2.44)]{kil}.

From the definition \Cref{jfp}, it follows that
\begin{align}
Q_n^{\alpha,\beta}(x) = \left(\frac{1+x}{2}\right)^{\alpha}Q_n^{0,\beta}(x), \label{Qn}
\end{align}
a relation we will frequently exploit. For brevity, we denote by $\mathbf{Q}^{\alpha,\beta}$ the quasimatrix of JFPs:
\begin{align*}
\mathbf{Q}^{\alpha,\beta} = \left(Q_{0}^{\alpha,\beta}(x), Q_{1}^{\alpha,\beta}(x), Q_{2}^{\alpha,\beta}(x), \ldots \right).
\end{align*}

Consider a function $u(x)$ that admits an expansion in the basis $\{Q_n^{\alpha,\beta}(x)\}_{n=0}^{\infty}$, i.e.,
\begin{align}
u(x) = \sum_{n=0}^{\infty} \hat{u}_n  Q_n^{\alpha,\beta}(x). \label{f}
\end{align}
The central objective of this paper is to construct an infinite-dimensional matrix $\mS$ such that
\begin{align}
\mathcal{I}^{\mu} u(x) = \mathbf{Q}^{\alpha,\beta} \mS \hat{u}, \label{matrep}
\end{align}
where $\hat{u} = \left(\hat{u}_0,\hat{u}_1, \hat{u}_2,\ldots\right)^{\top}$ is the coefficient vector. In practice, we compute finite-dimensional truncations of $\mS$.

Our work is not the first attempt to develop spectral approximations to the FIO. When $\mu$ is rational, a matrix approximation can be constructed based on its action on a direct sum space of weighted Jacobi polynomials \cite{hal}. For instance, when $\mu = 1/2$, the sum-space basis comprises the Legendre polynomials $\{P_n(x)\}_{n=0}^{\infty}$ and the $\sqrt{1+x}$-weighted Chebyshev polynomials of the second kind $\{U_n(x)\}_{n=0}^{\infty}$. Interleaving them yields
\begin{align}
\left( P_0(x) \big| \sqrt{1+x}U_0(x) \big| P_1(x) \big| \sqrt{1+x}U_1(x) \big| P_2(x) \big| \sqrt{1+x}U_2(x) \big| \cdots \right). \label{dsbasis}
\end{align}
This collection forms a frame in the infinite-dimensional setting and becomes a basis upon finite truncation. However, as noted in \cite{pu}, the elements in this basis are not orthogonal. Specifically, the set $\{\sqrt{1+x}U_n(x)\}_{n=0}^{\infty}$ is orthogonal with respect to the weight function $\sqrt{(1-x)/(1+x)}$, that is,
\begin{align*}
\int_{-1}^{1} \sqrt{1+x}U_{m}(x)\sqrt{1+x}U_{n}(x)\sqrt{\frac{1-x}{1+x}} \,\mathrm{d}x = \frac{\pi}{2}\delta_{mn},
\end{align*}
whereas the Legendre polynomials $\{P_n(x)\}_{n=0}^{\infty}$ are orthogonal with respect to the constant weight function $1$. Consequently, the sum-space basis lacks a unified weight and cannot be orthogonal. We further note that, in the special case $\alpha = 0$, the JFPs in \Cref{jfp} may equivalently be constructed by orthogonalizing a sum-space basis, such as the one given in \Cref{dsbasis}. The matrix approximation to the FIO in such a sum-space basis is banded \cite[\S 2.4]{hal}, but the lack of orthogonality in the basis can result in poor conditioning of its finite truncations. In solving FIEs or FDEs using such matrices, the resulting coefficient vectors may exhibit extremely large entries, requiring extended-precision arithmetic \cite[Example 3]{pu}. In such cases, the sum-space basis offers no advantage over the monomial basis $\{(1+x)^{n\mu}\}_{n=0}^{\infty}$ in terms of conditioning or coefficient magnitudes.

Pu and Fasondini propose in \cite[\S 5]{pu} two algorithms for constructing the infinite-dimensional matrix approximation to the FIO using JFPs as basis functions, i.e., $\mS$ in \Cref{matrep}. However, both algorithms rely on pseudo-stabilization techniques that require extended-precision arithmetic, leading to prohibitively high computational and storage costs. The first algorithm has a computational complexity of $\mathcal{O}(N^4 \log N \log \log N)$, where $N$ denotes the truncation size of $\mS$. The second algorithm achieves a lower complexity of $\mathcal{O}(N^3 \log N \log \log N/\beta)$, but it applies only when $1/\beta$, $\gamma - 1/\beta$, and $\alpha + 1/\beta - 1 - \gamma$ are all integers. Here, $\gamma$ is the second parameter in the Jacobi polynomial $P^{(\rotatebox{45}{\scalebox{0.4}{$\square$}}, \gamma)}(x)$, which replaces $T_n$ in \Cref{jfp}. This constraint implies that the second algorithm necessitates working with Jacobi polynomials of general fractional parameters, rather than the simpler Chebyshev or Legendre polynomials.

In this paper, we propose a new approach for constructing $\mS$ that is both fast and stable. Our method dispenses with extended-precision arithmetic, and achieves optimal complexity $\mathcal{O}(N^2)$ with negligible memory overhead. Importantly, the construction relies solely on Chebyshev polynomials, eliminating the need to handle Jacobi polynomials with general fractional parameters.

The availability of $\mS$ at such low computational cost can significantly benefit a wide range of applications. Most notably, as demonstrated in \cite{pu}, $\mS$ enables the JFP spectral method for solving fractional integral equations (FIEs) of the general form
\begin{align}
a_0(x) u(x) + a_1(x) \mathcal{I}^{\mu_1}\left[b_1(\rotatebox{45}{\scalebox{0.55}{$\square$}}) u\right](x) + \cdots + a_{\ell}(x) \mathcal{I}^{\mu_{\ell}}\left[b_{\ell}(\rotatebox{45}{\scalebox{0.55}{$\square$}}) u\right](x) = f(x), \label{fie}
\end{align}
where $\mu_1 > 0$, and $\mu_j = (p_j/q_j)\mu_1$ for $j = 2, \ldots, \ell$, with $p_j, q_j \in \mathbb{N}^+$ and $p_j/q_j$ irreducible. Let $q$ be the least common multiple of $\{q_j\}_{j=1}^{\ell}$ and $\tilde{\beta} = \mu_1/q$. It is also assumed that there exist $\tilde{\alpha}$ and $\check{\alpha}$ so that $a_j(x)b_j(x)\in \{Q_n^{\tilde{\alpha}, \tilde{\beta}}(x)\}_{n=0}^{\infty}$ for $j = 0, \ldots, \ell$, where $b_0(x) = 1$, and that $f(x)\in \{Q_n^{\check{\alpha}, \tilde{\beta}}(x)\}_{n=0}^{\infty}$.

Moreover, this JFP-based spectral method can also be applied to fractional differential equations (FDEs) that can be recast into the form of \Cref{fie} via integral reformulation techniques; see \cite{gre,du,pu} for details.

Apart from the JFP spectral method, there exist four other spectral methods for solving FIEs and FDEs: the collocation method based on polyfractonomials (PFC) \cite{zay}, the Petrov--Galerkin method using generalized Jacobi functions (GJFPG) \cite{che}, the sum-space method (SS) \cite{hal}, and the Petrov--Galerkin method based on generalized log-orthogonal functions (GLOFPG) \cite{che2,che5}.

The first three methods all use certain variants of the weighted Jacobi polynomials as (part of) the basis. The PFC and GJFPG methods approximate the algebraic singularities in the solution using polynomials, or polynomials in fractional powers that differ from those of the exact solution. Consequently, they generally fail to achieve spectral convergence, and the convergence is typically only algebraic \cite{che,pu,zay}.
An additional downside of the PFC method is the system it leads to is dense, whereas the system due to the GJFPG method may be sparse for certain select problems. The SS method addresses these issues by employing bases that span the direct sum of suitably weighted ultraspherical and Jacobi polynomial spaces---it usually converges exponentially and the resulting system is banded or lower-banded\footnote{Lower-banded matrices are sometimes referred to as $m$-Hessenberg matrices, where $m$ indicates the lower bandwidth.} for FIEs and FDEs with constant and variable coefficients respectively \cite{hal}. As noted above, the linear system, however, may be ill-conditioned \cite{pu}. Moreover, the SS method can only handle FIEs and FDEs of rational order, i.e., $\mu = p/q$, where $p$ and $q$ are positive integers. It resorts to $2q$ different weighted orthogonal polynomial bases. When $q$ is not small, the method becomes increasingly unwieldy to implement. Finally, we note that the PFC, GJFPG, and SS methods all rely on the use of Jacobi polynomials. For Jacobi polynomials with general integer or fractional parameters, the transforms between values and coefficients are much less efficient than those for Chebyshev or Legendre polynomials\footnote{Even for Legendre polynomials, our experience with the existing algorithms for transforms between values and coefficients is far from satisfactory. The algorithms are galactic due to the powers of $\log n$ in the asymptotic complexity and the huge hidden constant in the big-Oh notation.}. 

The GLOFPG method \cite{che2} represents the solution to FDEs using generalized log-orthogonal functions (GLOF). These functions are capable of approximating weak singularities over a much broader class than JFPs, thereby allowing the GLOFPG method to handle problems beyond those conforming to \Cref{fie} or reducible to it. Despite this generality, the method suffers from several drawbacks: (1) The construction of the associated linear system is of complexity $\mathcal{O}(N^4)$, as the system matrix is dense and the evaluation of each entry requires $\mathcal{O}(N^2)$ flops. (2) The condition number of the system matrix deteriorates rapidly with increasing system size. For matrices of just a few hundred in size, the system can become so ill-conditioned that the computed solution may contain no reliable digits. (3) The evaluation of high-degree GLOFs is prone to underflow or overflow, further limiting the maximum feasible truncation size. 


The JFP spectral method, empowered by the fast and stable construction of $\mS$ proposed in this paper, significantly outperforms the PFC, GJFPG, SS, and GLOFPG methods in solving FIEs and FDEs. A summary of the key features and limitations of these methods is provided in \Cref{tab:compare}.

\begin{table}[t!]
  \centering
  \renewcommand{\arraystretch}{1.4}
  \setlength{\tabcolsep}{8pt}
  \caption{Comparison of three mainstream spectral methods for solving general variable-coefficient FIEs and FDEs. For the SS method, the structure and computational cost of the resulting linear system depend on whether the coefficients are constant or variable.}\label{tab:compare}
  \begin{tabular}{c c c c c}
    \toprule
    \multicolumn{2}{c}{method} & SS & GLOFPG & JFP \\
    \midrule
    \multicolumn{2}{c}{$\mu$} & rational & no restriction & no restriction \\
    \midrule
    \multirow{2}{*}{structure} & constant & strictly banded & \multirow{2}{*}{full} & \multirow{2}{*}{lower-banded} \\
                                  & variable & lower-banded &  &  \\
    \midrule
    \multirow{2}{*}{construction} & constant & $\mathcal{O}(N)$ & \multirow{2}{*}{$\mathcal{O}(N^4)$} & \multirow{2}{*}{$\mathcal{O}(N^2)$} \\
                                  & variable & $\mathcal{O}(N^2)$ &  &  \\
    \midrule
    \multirow{2}{*}{solve} & constant & $\mathcal{O}(N)$ & \multirow{2}{*}{$\mathcal{O}(N^3)$} & \multirow{2}{*}{$\mathcal{O}(N^2)$} \\
                                  & variable & $\mathcal{O}(N^2)$ &  &  \\                              
    \midrule
    \multicolumn{2}{c}{conditioning} & ill & ill & good \\
    \bottomrule
  \end{tabular}
\end{table}

The applications of the spectral approximation to the FIO extend well beyond the solution of boundary value problems for FIEs and FDEs. In some cases, the solutions to certain FIEs are Mittag--Leffler functions, indicating that the JFP spectral method may serve as an effective tool for computing these special functions \cite{pu}. Moreover, the matrix approximations developed here can be applied to the computation of eigenvalues and pseudospectra of FIOs and fractional differential operators (FDOs). In the context of one-dimensional FDE initial value problems or time-dependent fractional partial differential equations, these spectral approximations enable time-fractional derivatives to be integrated with spectral accuracy via the spectral deferred correction method \cite{dut}, while using only low-order fractional time-stepping schemes.

Throughout this paper, the two-parameter Mittag--Leffler function is denoted by $E_{\sigma, \tau}$, following standard convention. We shall say that a matrix $A$ has bandwidths $(\xi_l, \xi_u)$ if the entries of $A$ satisfy $A_{ij}=0$ for $i-j>\xi_l$ and $j-i>\xi_u$. $\Re( \rotatebox{45}{\scalebox{0.5}{$\square$}} )$ and $\Im( \rotatebox{45}{\scalebox{0.5}{$\square$}} )$ denote the real and imaginary parts respectively.

The remainder of the paper is organized as follows. In \Cref{sec:con}, we present a fast and stable method for constructing the spectral approximation to the FIO. Various applications discussed above are demonstrated in \Cref{sec:examples}. We conclude with a brief outlook on future work.

\section{Constructing the matrix}\label{sec:con}
For $u(x)$ given in \Cref{f}, the action of the FIO on $u(x)$ amounts to that on $Q_n^{\alpha,\beta}(x)$ for each $n$. In the language of quasimatrices, this is 
\begin{align}
\mathcal{I}^{\mu}u(x) = \mathcal{I}^{\mu} \mathbf{Q}^{\alpha,\beta}\hat{u} = \left( 
    \vphantom{\begin{aligned}
      & \\ 
      & \\ 
      &
  \end{aligned}}
  \right.
  \begin{aligned}
    \mathcal{I}^{\mu}Q_0^{\alpha,\beta}(x)\Bigg|  \mI^{\mu}Q_1^{\alpha,\beta}(x) \Bigg| \mI^{\mu}Q_2^{\alpha,\beta}(x) \Bigg| \cdots  
  \end{aligned}
  \left.
  \vphantom{\begin{aligned}
      & \\ 
      & \\ 
      &
  \end{aligned}}
  \right)
\hat{u}. \label{quasimat}
\end{align}
Consider the $n$th column in this quasimatrix, that is,
\begin{align*}
\mathcal{I}^{\mu}Q_n^{\alpha,\beta}(x) &= \frac{2^{-\alpha}}{\Gamma(\mu)} \int_{-1}^{x} \frac{(1+t)^{\alpha} }{(x-t)^{1-\mu}}T_n\left(2\left(\frac{1+t}{2}\right)^{\beta}-1\right) \md t.
\end{align*}
Let us make a change of variable $t = x - (1+x)(1+s)/2$ to have
\begin{align*}
\mathcal{I}^{\mu}Q_n^{\alpha,\beta}(x) = \frac{2^{-\alpha}}{\Gamma(\mu)} \left(\frac{1+x}{2}\right)^{\alpha +\mu}  \varphi_n(x), 
\end{align*}
where
\begin{align}
\varphi_n(x) =  \int_{-1}^{1} \frac{(1-s)^{\alpha}}{(1+s)^{1-\mu}} T_n\left(2\left(\frac{(1+x)(1-s)}{4}\right)^{\beta} -1\right) \md  s. \label{phi}
\end{align}
The effect of this change of variable is twofold---$x$ and $t$ are now decoupled, and the upper limit of integration is no longer variable. A key observation to make is that $\varphi_n(x)$ is a polynomial in $(1+x)^{\beta}$ of degree at most $n$. Thus, it can be written as
\begin{align*}
\varphi_n(x) = \sum_{m=0}^{n} R_{mn} Q_m^{0,\beta}(x),
\end{align*}
where $\{R_{mn}\}_{m=0}^n$ are the coefficients. As shown below by \Cref{thm:M}, the pre-multiplication of $\left(\frac{1+x}{2}\right)^{\mu}$ can be represented by the infinite-dimensional multiplication matrix $\mM$ such that 
\begin{align}
\left(\frac{1+x}{2}\right)^{\mu} \mathbf{Q}^{0,\beta} \hat{g} = \mathbf{Q}^{0,\beta} \mM \hat{g}, \label{multiplication}
\end{align}
for an infinite vector $\hat{g}$. It then follows from \Cref{quasimat}, \Cref{Qn}, and the last two equations that 
\begin{align*}
\mathcal{I}^{\mu}u(x) = \frac{2^{-\alpha}}{\Gamma(\mu)} \mathbf{Q}^{\alpha,\beta}\mM\mR\hat{u},
\end{align*}
where $\mathcal{R}$ is an infinite-dimensional upper triangular matrix with entries $R_{mn}$ for $m, n = 0, 1, \ldots$. Thus, 
\begin{align}
\mS = \frac{2^{-\alpha}}{\Gamma(\mu)} \mM \mathcal{R} \label{S}
\end{align}
is the infinite-dimensional matrix approximation to the FIO defined in \Cref{matrep}. Hence, our task now boils down to the construction of $\mR$ and $\mM$.

In the rest of this section, we shall denote by $\mathcal{R}_n$ the $n$th column of $\mathcal{R}$ with the index $n$ starts from 0. Now we start off by showing the recurrence relation satisfied by $\varphi_n(x)$. 

\subsection{Recurrence relation satisfied by $\varphi_n(x)$}\label{sec:rec}
We begin by recalling some of the basic properties of the Chebyshev polynomials, which can be found in many standard texts, e.g., \cite{sze}.
\begin{lemma}
For Chebyshev polynomials of the first kind $T_n(x)$ and the second kind $U_n(x)$,
\begin{subequations}
\begin{align}
(n-1)\frac{\md}{\md x}T_{n+1}(x) &= 2(n^2-1)T_{n}(x) + (n+1)\frac{\md}{\md x}T_{n-1}(x), \label{DTT}\\
T_{n+1}(x) &= 2xT_{n}(x)-T_{n-1}(x), \label{Trec}\\
\frac{\md}{\md x}T_{n}(x) &= nU_{n-1}(x), \label{DTU}\\
U_{n+1}(x) &= 2xU_n(x)-U_{n-1}(x), \label{Urec}\\
T_m(x)T_n(x) &=\frac{1}{2}\left(T_{m+n}(x)+T_{|m-n|}(x)\right), \label{TT}
\end{align}
where $n \geq 1$ and $m \geq 0$.
\end{subequations}
\end{lemma}
We now present the main result of this section, establishing that $\varphi_{n+1}(x)$ satisfies a three-term recurrence relation in the form of a differential equation.

\begin{theorem}(recurrence relation) For $n\geq 2$,
\begin{align}
\left(\frac{1+x}{n+1}\frac{\md}{\md x} - \beta\right) \varphi_{n+1} (x) = 2\beta \varphi_{n} (x) + \left(\beta+\frac{1+x}{n-1}\frac{\md}{\md x}\right) \varphi_{n-1} (x). \label{rec}
\end{align} 
\end{theorem}
\begin{proof}
To make the derivation uncluttered, let 
\begin{align}
\tilde{y} = 2\left(\frac{(1+x)(1-s)}{4}\right)^{\beta} - 1. \label{xymap}
\end{align}
By \Cref{DTT}, we have
\begin{align}
\frac{\md}{\md x} \varphi_{n+1} (x) & = \int_{-1}^{1} \frac{(1-s)^{\alpha}}{(1+s)^{1-\mu}} \frac{\md }{\md x} T_{n+1}\left(\tilde{y}\right) \md s \nonumber \\
&= 2^{2-2\beta} \beta (n+1) (1+x)^{\beta-1} \int_{-1}^{1} \frac{(1-s)^{\alpha}}{(1+s)^{1-\mu}} (1-s)^{\beta}T_{n}\left(\tilde{y}\right) \md s \label{interm} \\
&\quad + \frac{n+1}{n-1} \int_{-1}^{1} \frac{(1-s)^{\alpha}}{(1+s)^{1-\mu}} \frac{\md }{\md x} T_{n-1}\left(\tilde{y}\right) \md s \nonumber,
\end{align}
where the last term
\begin{align}
\frac{n+1}{n-1}\int_{-1}^{1} \frac{(1-s)^{\alpha}}{(1+s)^{1-\mu}} \frac{\md }{\md x} T_{n-1}\left(\tilde{y}\right) \md s = \frac{n+1}{n-1}\frac{\md}{\md x} \varphi_{n-1}(x). \label{term2}
\end{align}
By \Cref{xymap}, the first term on the right-hand side of the second equality in \Cref{interm}
\begin{align}
  \begin{aligned}
  2^{2-2\beta}\beta(n+1) (1+x)^{\beta-1} &\int_{-1}^{1} \frac{(1-s)^{\alpha}}{(1+s)^{1-\mu}} (1-s)^{\beta}T_{n}\left(\tilde{y}\right) \mathrm{d} s \\
  &= \frac{2\beta(n+1)}{1+x} \int_{-1}^{1} \frac{(1-s)^{\alpha}}{(1+s)^{1-\mu}} (\tilde{y}+1)T_{n}\left(\tilde{y}\right) \mathrm{d} s. \label{term1}
  \end{aligned}
  \end{align}
Using \Cref{Trec} gives
\begin{align}
  \begin{aligned}
  \int_{-1}^{1} \frac{(1-s)^{\alpha}}{(1+s)^{1-\mu}} \tilde{y}T_{n}\left(\tilde{y}\right) \mathrm{d} s &= \frac{1}{2}\int_{-1}^{1} \frac{(1-s)^{\alpha}}{(1+s)^{1-\mu}}\left(T_{n+1}(\tilde{y})+T_{n-1}(\tilde{y})\right)\mathrm{d} s \\
  &= \frac{1}{2}\left(\varphi_{n+1} (x)+\varphi_{n-1} (x)\right). \label{term1b}
  \end{aligned}
  \end{align}
Substituting \Cref{term2,term1,term1b} back to \Cref{interm} gives \Cref{rec}.
\end{proof}

The recurrence relation \Cref{rec} shows that $\varphi_{n+1}(x)$ can be obtained once $\varphi_{n-1} (x)$ and $\varphi_n(x)$ are available for $n \geq 2$. We now turn to the calculation of $\varphi_0(x)$, $\varphi_1(x)$, and $\varphi_2(x)$, which initiates the recursion. The following lemma prepares the key value $h_n$, which shall be frequently used.
\begin{lemma}
\begin{align*}
h_n = \int_{-1}^{1}\frac{(1-s)^{\alpha + n \beta}}{(1+s)^{1-\mu}}\md s = 2^{\mu+\alpha+n\beta} \mathrm{B}(\mu,1+\alpha+n\beta),
\end{align*}
where $\mathrm{B}(\cdot)$ is the beta function.
\end{lemma}
\begin{proof}
By change of variable $s = 2\tilde{s}-1$, we have
\begin{align*}
h_n = 2^{\mu+\alpha+n\beta} \int_{0}^{1} (1-\tilde{s})^{\alpha+n\beta}\tilde{s}^{\mu-1} \md \tilde{s} = 2^{\mu+\alpha+n\beta} \mathrm{B}(\mu,1+\alpha+n\beta).
\end{align*}
\end{proof}
\begin{theorem}
The first three $\varphi_n(x)$ can be expressed in terms of $h_n$ and $\{Q_n^{0,\beta}(x)\}_{n=0}^{2}$ as follows.
\begin{subequations}
\begin{align}
\varphi_{0}(x) &= h_0,\label{phi0} \\
\varphi_{1}(x) &= \frac{h_1}{2^{\beta}} - h_0 + \frac{h_1}{2^{\beta}} Q_1^{0,\beta}(x), \label{phi1}\\
\varphi_{2}(x) &= \frac{3h_2}{2^{2\beta}} -\frac{h_1}{2^{\beta-2}} + h_0 + 4\left(\frac{h_2}{2^{2\beta}} - \frac{h_1}{2^{\beta}}\right) Q_1^{0,\beta}(x) + \frac{h_2}{2^{2\beta}} Q_2^{0,\beta}(x).\label{phi2}
\end{align}\label{phi_0}
\end{subequations}
\end{theorem}
\begin{proof}
It is straightforward to see \Cref{phi0}. By definition,
\begin{align*}
\varphi_1(x) &=  \int_{-1}^{1} \frac{(1-s)^{\alpha}}{(1+s)^{1-\mu}} \left(2\left(\frac{(1+x)(1-s)}{4}\right)^{\beta} -1\right) \md s\\
& = 2^{1-2\beta} (1+x)^{\beta} \int_{-1}^{1} \frac{(1-s)^{\alpha+\beta}}{(1+s)^{1-\mu}} \md s -  \int_{-1}^{1} \frac{(1-s)^{\alpha}}{(1+s)^{1-\mu}} \md s\\
& = 2^{1-2\beta}h_1  (1+x)^{\beta} - h_0,
\end{align*}
which gives \Cref{phi1}. We omit the proof of $\varphi_2(x)$, for it is similar to that of $\varphi_1(x)$.
\end{proof}

With \Cref{phi1,phi2}, we can recurse for $\varphi_n(x)$ for $n\geq 3$ using \Cref{rec}. Yet one thing is still missing---we need a boundary condition to anchor the solution $\varphi_{n+1}(x)$, as \Cref{rec} is a first-order ODE of $\varphi_{n+1}(x)$. To this end, we take the value of $\varphi_{n+1}(x)$ at $x=1$, i.e., 
\begin{align}
\varphi_{n+1} (1) = \int_{-1}^{1} \frac{(1-s)^{\alpha}}{(1+s)^{1-\mu}} T_{n+1}\left(2\left(\frac{1-s}{2}\right)^{\beta}-1\right) \md s. \label{inival}
\end{align}
as the boundary condition\footnote{For \Cref{rec}, the conditions of the Picard--Lindel\"{o}f theorem are not satisfied when the left boundary condition is prescribed; therefore, the existence and uniqueness of a solution are not guaranteed.}. Noting that the integrand exhibits weak singularities at both endpoints, we evaluate \Cref{inival} using the double-exponential method \cite{tak}, which transforms the original integral into a doubly infinite one with a rapidly decaying integrand.

Our extensive experiments show that using the trapezoidal rule with either $8n$ or $80$ equispaced points in $[-4, 4]$, whichever is greater, suffices to evaluate \Cref{inival} to machine precision, at a cost of $\mathcal{O}(n)$ or less.

\subsection{Recursing for $\mR$}\label{sec:R} 
Now we have all the ingredients to generate $\varphi_n(x)$ recursively following \Cref{rec}. The remaining task is to solve the ODE for the $\mathbf{Q}^{0,\beta}$ coefficients of $\varphi_n(x)$, i.e., $R_{mn}$ for $m = 0, 1, \ldots, n$. The relationship \Cref{jfp} between the JFPs and the Chebyshev polynomials of the first kind motivates the solution of \Cref{rec} using a variant of the ultraspherical spectral method tailored to the basis $\mathbf{Q}^{0,\beta}$.

We first consider the matrix approximation to the weighted differential operator, which appears on both sides of \Cref{rec}. For convenience, define
\begin{align}
y = 2\left(\frac{x+1}{2}\right)^{\beta} - 1, \label{cov}
\end{align}
and let $\mathbf{T} = \left(1, T_1(y), T_2(y), \dots \right)$ and $\mathbf{U} = \left(1, U_1(y), U_2(y), \dots \right)$ denote the quasimatrices of Chebyshev polynomials of the first and second kinds, respectively, expressed in terms of the variable $y$. Note that $\mathbf{T} = \mathbf{Q}^{0,\beta}$. 

\begin{lemma}\label{lem:wdo}
For the weighted differential operator $\displaystyle (1+x)\frac{\md}{\md x}$, we have
\begin{align*}
  (1+x)\frac{\md}{\md x}\mathbf{T} = \beta \mathbf{U}
  \underbrace{
    \begin{pmatrix*}[r]
        0 & 1 & 1 & \\
          & \frac{1}{2} & 2 & \frac{3}{2}\\
          &   & 1 & 3 & 2 \\
          & & & \ddots& \ddots& \ddots
    \end{pmatrix*}
  }_{\mathclap{\text{\large $\mathcal{D}$}}},
\end{align*}
where $\mD$ is a banded matrix of infinite dimensions with bandwidths $(0,2)$.
\end{lemma}
\begin{proof}
By \Cref{DTU}, we have
\begin{align*}
(1+x)\frac{\md}{\md x }T_n\left(y \right) &= (1+x)\left(2^{1-\beta}n\beta (1+x)^{\beta-1} U_{n-1}\left(y \right)\right)\\
& = n \beta \left(y+1\right) U_{n-1}\left(y \right)\\
& = \frac{n\beta}{2}\left(U_{n-2}\left(y \right) + 2U_{n-1}\left(y \right)+ U_{n}\left(y \right)\right),
\end{align*}
where in the last equality we have used \Cref{Urec}.
\end{proof}

What is also required is the matrix approximation to the conversion operator that maps a series of Chebyshev polynomials of the first kind to that of the second, both in $y$. We omit the proof, since this conversion operator is identical to the one used in the standard ultraspherical spectral method \cite{olv}.
\begin{lemma}\label{lem:conversion}
For $\mathbf{T}$ and $\mathbf{U}$,
\begin{align*}
\mathbf{T} = \mathbf{U}
\underbrace{\begin{pmatrix}
  1 & &-\frac{1}{2}\\
  & \frac{1}{2} & & -\frac{1}{2}\\
  & &\frac{1}{2} & & -\frac{1}{2}\\
  & & &\ddots& &\ddots\\
\end{pmatrix}}_{\mathclap{\text{\large $\mathcal{C}$}}},
\end{align*}
where $\mC$ is the infinite-dimensional conversion matrix with bandwidths $(0, 2)$.
\end{lemma}

To represent the Dirichlet boundary condition at the right endpoint, we need an infinite row vector with each entry being the value of $Q_n^{0,\beta}(1)$.
\begin{lemma}\label{lem:dbc}
The action of the Dirichlet boundary condition on $\mathbf{T}$ can be represented by 
\begin{align}
\mB = \left(1,1,1,\cdots\right). \label{bc}
\end{align}
\end{lemma}
\begin{proof}
This follows from $Q_n^{0,\beta}(1) =  T_n\left(1\right) = 1$.
\end{proof}

Following \Cref{lem:wdo,lem:conversion,lem:dbc}, we can formulate an infinite-dimensional linear system that represents \Cref{rec}. However, since it is known \textit{a priori} that $\mR_j$ contains exactly $n+1$ nonzero entries, we can instead work with a finite $(n+2)\times(n+2)$ system to solve for the nonzero entries of $\mR_{n+1}$, without resorting to any adaptive procedure to determine an optimal truncation size for the infinite system. Let $D$ and $C$ denote the $(n+1)\times (n+2)$ truncations of $\mD$ and $\mC$, respectively, and let $R_j$ and $B$ denote the $(n+2)$-term truncations of $\mathcal{R}_j$ for $j \in \{n-1,n,n+1\}$ and $\mB$, respectively. Then, \Cref{rec} can be discretized as
\begin{align}
\begin{pmatrix*}[c]
B \\[1mm]
\displaystyle \frac{1}{n+1}D - C
\end{pmatrix*}R_{n+1}=
\begin{pmatrix*}[c]
\varphi_{n+1}(1) \\[4mm]
\displaystyle G
\end{pmatrix*}, \label{system}
\end{align}
where
\begin{align*}
G = 2 CR_{n} + \left( C + \frac{1}{n-1}D\right)R_{n-1}.
\end{align*}
As in the standard ultraspherical spectral method, \Cref{system} is an almost-banded system with bandwidths $(1,1)$, which can be solved in linear complexity. However, a recombined basis that satisfies the Dirichlet boundary condition can be employed to obtain a strictly banded system of bandwidths $(1,2)$, for which a standard banded solver can be used for a significant speed boost. For details, see, for example, \cite{qin}.

For the nonzero entries in the first $N+1$ columns of $\mR$, the total cost of the recursion is $\mO(N^2)$ flops. 

\subsection{Adaptive construction of $\mS$}\label{sec:adapt}
Up to this point, we are just one step away from $\mS$---constructing $\mM$ in \Cref{multiplication} to represent the action of pre-multiplication of $\left((1+x)/2\right)^{\mu}$. Since $\mu = k\beta$ for some positive integer $k$, $\left((1+x)/2\right)^{\mu}$ is the $k$th power of $\left((1+x)/2\right)^{\beta}$. It therefore can be written as a linear combination of $Q_n^{0,\beta}(x)$, i.e.,
\begin{align}
\left(\frac{1+x}{2}\right)^{\mu} = \sum_{j=0}^{k} c_j Q_j^{0,\beta}(x),\label{mulexp}
\end{align}
or, equivalently,
\begin{align*}
\left(\frac{1+y}{2}\right)^k = \sum_{j=0}^{k} c_j T_j(y),
\end{align*}
where $y$ is defined in \Cref{cov}. The coefficients $\{c_j\}_{j=0}^k$ can be determined using discrete cosine transform (DCT) or FFT at a cost of $\mO(k\log k)$ \cite{atap}. Alternatively, we can calculate $\{c_j\}_{j=0}^k$ using \Cref{TT} and the fact that $(1+y)/2 = (T_0(y)+T_1(y))/2$; the cost is $\mO(k)$. 

The follow lemma shows that the infinite-dimensional multiplication matrix that represents multiplication of an infinite $Q^{0, \beta}_j(x)$ series by another one is Toeplitz plus Hankel up to a rank one perturbation, identical to that in the standard ultraspherical spectral method.
\begin{lemma}\label{thm:M}
Multiplying the infinite series $\sum_{j=0}^{\infty} c_j Q_j^{0,\beta}(x)$ by another infinite $Q^{0, \beta}_j(x)$ series can be represented as
\begin{align}
\mM = \frac{1}{2}\left[
\begin{pmatrix*}[r]
2c_0 & c_1 & c_2 & c_3 & \cdots\\
c_1 & 2c_0 & c_1 & c_2 & \ddots\\
c_2 & c_1 & 2c_0 & c_1 & \ddots\\
c_3 & c_2 & c_1 & 2c_0 & \ddots\\
\vdots & \ddots &\ddots & \ddots & \ddots
\end{pmatrix*} +
\begin{pmatrix*}
0 & 0 & 0 & 0 & \cdots\\
c_1 & c_2 & c_3 & c_4 & \ddots\\
c_2 & c_3 & c_4 & c_5 & \ddots\\
c_3 & c_4 & c_5 & c_6 & \ddots\\
\vdots & \ddots &\ddots & \ddots & \ddots
\end{pmatrix*}
\right].
\label{M}
\end{align}
\end{lemma}
\begin{proof}
By the same change of variable \Cref{cov} and \Cref{TT}, we have
\begin{align*}
Q_m^{0,\beta}(x) Q_n^{0,\beta}(x) &= T_m(y) T_n(y) = \frac{1}{2}\left(T_{m+n}(y)+T_{\left\lvert m-n\right\rvert}(y)\right)\\
 &= \frac{1}{2}\left(Q_{m+n}^{0,\beta}(x)+Q_{\left\lvert m-n\right\rvert }^{0,\beta}(x)\right).
\end{align*}
Thus, the same multiplication matrix as in \cite{olv} follows.
\end{proof}
For $((1+x)/2)^{\mu}$, whose expansion in terms of $Q_n^{0,\beta}(x)$ is a series of degree $k$, the matrix $\mM$ has bandwidths $(k, k)$. Since $\mR$ is upper triangular, the infinite-dimensional matrix $\mS$ has bandwidths $(k, \infty)$. Consequently, constructing the $(N+1)\times(N+1)$ truncation of $\mM$ requires $\mO(N)$ operations for $N \gg k$, and forming an $N\times N$ finite section of $\mS$ according to \Cref{S} incurs a computational cost of $\mO(N^2)$. The overall procedure for constructing a finite section of $\mS$ is summarized in \Cref{alg}.

\begin{algorithm}[t!]
\caption{Construction of the matrix approximation to a FIO}\label{alg}
\begin{algorithmic}
\STATE{Determine $\alpha$ and $\beta$ so that $\mu = k\beta$ for the smallest integer $k$.}
\STATE{Calculate the nonzero entries in $\mR_j$ for $j = 0, 1, 2$ following \Cref{phi}.}
\FOR{$n = 3$ to $N$}
    \STATE{Calculate the boundary condition \Cref{inival}.}
    \STATE{Solve ODE \Cref{system} to obtain the nonzero entries in $\mathcal{R}_n$.}
\ENDFOR
\STATE{Calculate $c_n$ in \Cref{mulexp} and construct the top-left $(N+1)\times(N+1)$ finite section of the multiplication matrix $\mM$.}
\STATE{Form the $(N+1)\times(N+1)$ approximation to the FIO following \Cref{S}.}
\end{algorithmic}
\end{algorithm}

In many of the applications discussed in \Cref{sec:intro}, such as solving FIEs of the form \Cref{fie}, we encounter a linear system whose coefficient matrix is constructed from finite sections of $\mS$ and is therefore lower-banded. Such systems can be efficiently solved using adaptive QR factorization in $\mO(N^2)$ flops. As in the standard ultraspherical spectral method, the residual can be evaluated by computing the norm of the nontrivial part of the right-hand-side vector. The iteration is terminated once the residual falls below a prescribed tolerance. The overall cost of this process is minimized when $\mS$ is constructed in an adaptive manner as well. Specifically, when the truncated system is enlarged to accommodate a longer solution, only the newly introduced rows and columns in $\mS$ and, in turn the corresponding nonzero entries in $\mM$ and $\mR$, need to be computed.

\section{Numerical examples}\label{sec:examples}
We illustrate the applicability of the spectral approximations to FIOs through several numerical examples.

\subsection{FIE boundary value problems}
\subsubsection{Fractional Abel integral equation}
\begin{figure}[t!]
\centering
\subfloat[]{\includegraphics[width=0.48\linewidth]{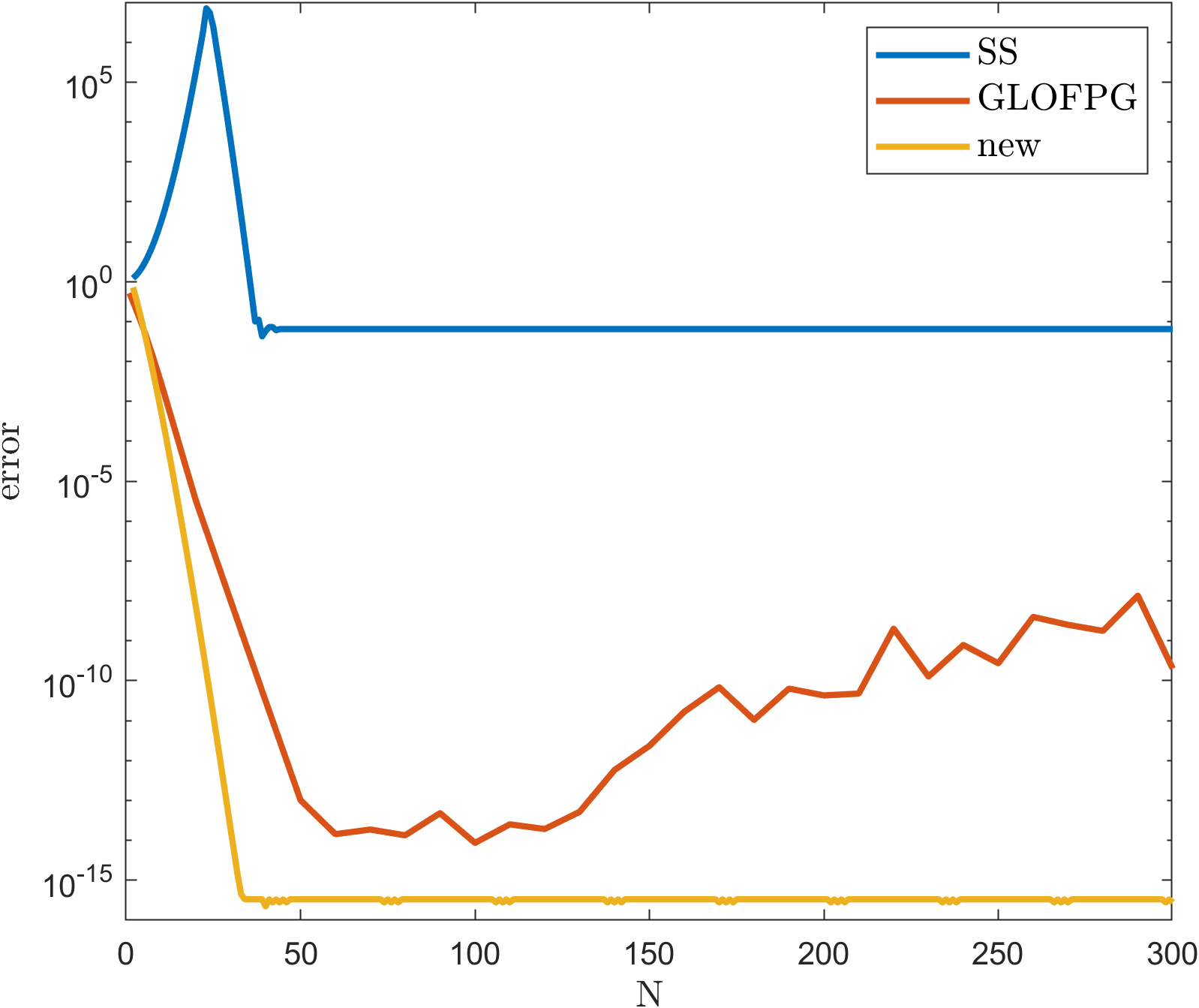}\label{fig:err}}
\hfill
\subfloat[]{\includegraphics[width=0.467\linewidth]{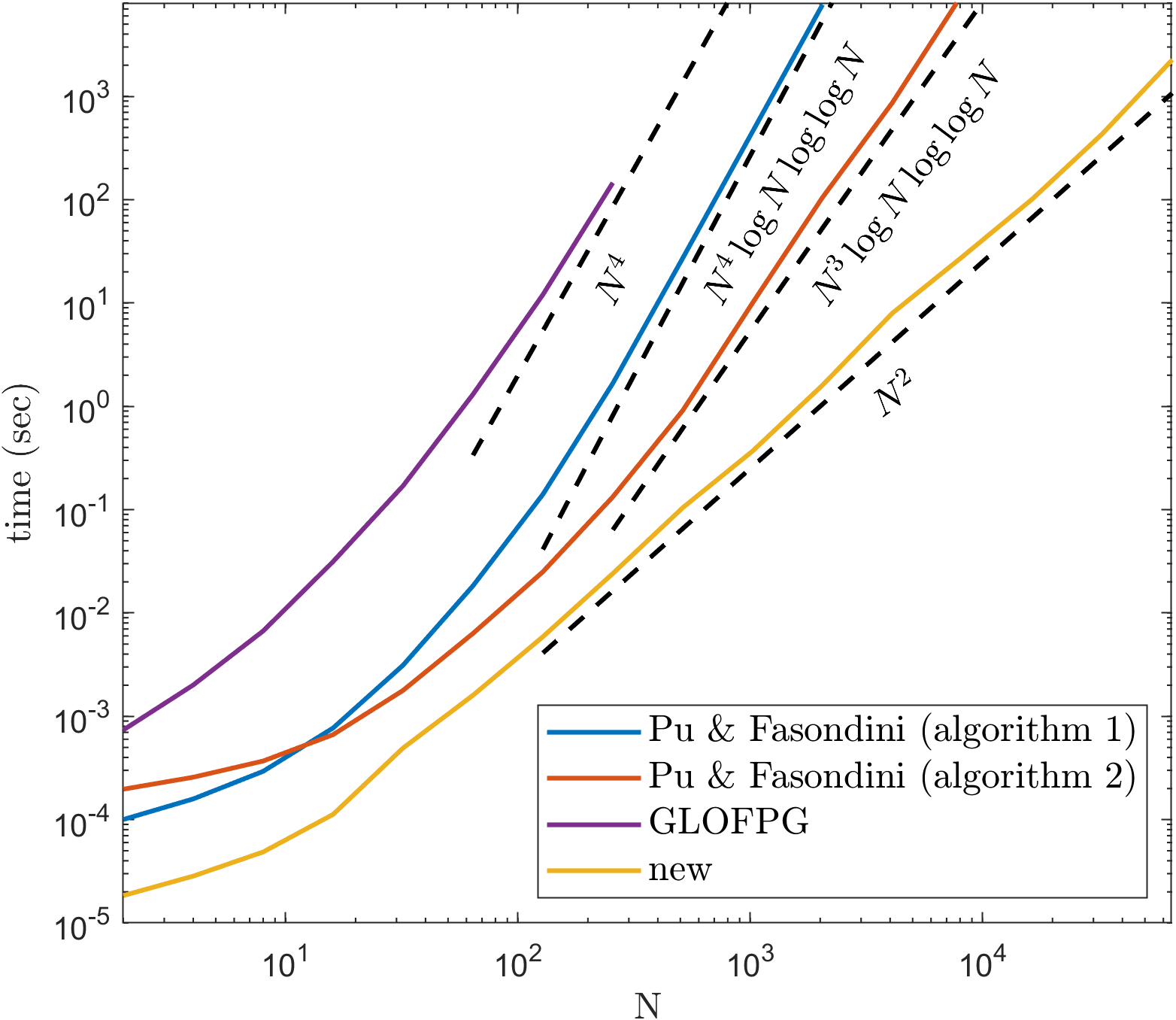}\label{fig:time}}
\caption{(A) Errors in the solutions to \Cref{abel} for $\lambda = 2$ obtained by the SS, JFP, and GLOFPG methods. For the GLOFPG method, we use the GLOF basis with parameters $\alpha = 0$, $\beta = 5$, and $\lambda = 0$, where $\alpha$, $\beta$, and $\lambda$ follow the notations in \cite[\S 3]{che2} and are not those used elsewhere in this paper. (B) Execution times for the GLOFPG method and for the JFP method constructed via different algorithms.}\label{fig:fie}
\end{figure}
First, we consider the second-kind Abel integral equation of fractional order 
\begin{align}
u(x) + \lambda^2 \mI^{1/2} u(x) = 1, ~~~ x \in[-1,1], \label{abel}
\end{align}
whose many applications can be found in \cite{pol}. The exact solution to \Cref{abel} can be written in terms of the Mittag-Leffler function  
\begin{align*}
u(x) = E_{\frac{1}{2},1}(-\lambda^2\sqrt{1+x}).
\end{align*}
The performance of the SS and JFP methods applied to \Cref{abel} is examined in detail in \cite{pu}, with particular attention to the large-$\lambda$ regime. What has been observed is as follows. For $\lambda = 1$, both methods work well. For the SS method, the condition number of the linear system grows as $\mO\left(\exp ({2\lambda^4})\right)$ as $\lambda$ increases. As a consequence, the largest coefficient of $u(x)$ grows in the same rate. For $\lambda = 2$, the largest coefficient of the solution is $\mathcal{O}(10^{14})$, so double-precision arithmetic can barely handle the computation; the severe cancellation errors limit the accuracy to no more than two digits. Beyond $\lambda = 2$, one has to resort to extended precision arithmetic for the SS method to produce any meaningful results. By contrast, in the JFP method the solution coefficients remain bounded by $1$ regardless of the value of $\lambda$, and the linear systems are much better conditioned with the condition number plateauing at approximately $1.85\lambda^2$ as $N \to \infty$; see \cite{pu}.

We solve \Cref{abel} for $\lambda = 2$ using the SS, JFP, and GLOFPG methods, and present the results in \Cref{fig:fie}. In \Cref{fig:err}, the error is plotted against the truncation size of the system. As expected, the SS method diverges initially before converging to an accuracy of $\mathcal{O}(10^{-2})$, after which it stagnates, whereas both the JFP and GLOFPG methods exhibit exponential convergence. Once convergence to machine precision is attained, the accuracy of the JFP method remains stable. This is attributed to the fact that the linear system arising from the JFP method is lower-banded \cite{che4}, which ensures the error does not bounce back. In contrast, the GLOFPG method, which leads to a dense linear system, fails to converge beyond approximately $\mathcal{O}(10^{-14})$, at which point the effects of ill-conditioning become dominant. By $n = 300$, nearly half of the significant digits have already been lost.

In \Cref{fig:time}, we present the total execution times for solving \Cref{abel}---including both system construction and solution---for the GLOFPG method ($\mathcal{O}(N^4)$) and the JFP spectral method, using three construction algorithms: the two from \cite{pu} with complexities $\mathcal{O}(N^3 \log N \log \log N)$ and $\mathcal{O}(N^4 \log N \log \log N)$ respectively, and our new algorithm with $\mathcal{O}(N^2)$. The GLOFPG method terminates at around $n = 300$ due to overflow when evaluating the GLOF basis. The SS method is excluded, as it is considered impractical for this problem.

\subsubsection{Integral equation of multiple fractional orders with variable coefficients}
\begin{figure}[t!]
\centering
\includegraphics[width=0.48\linewidth]{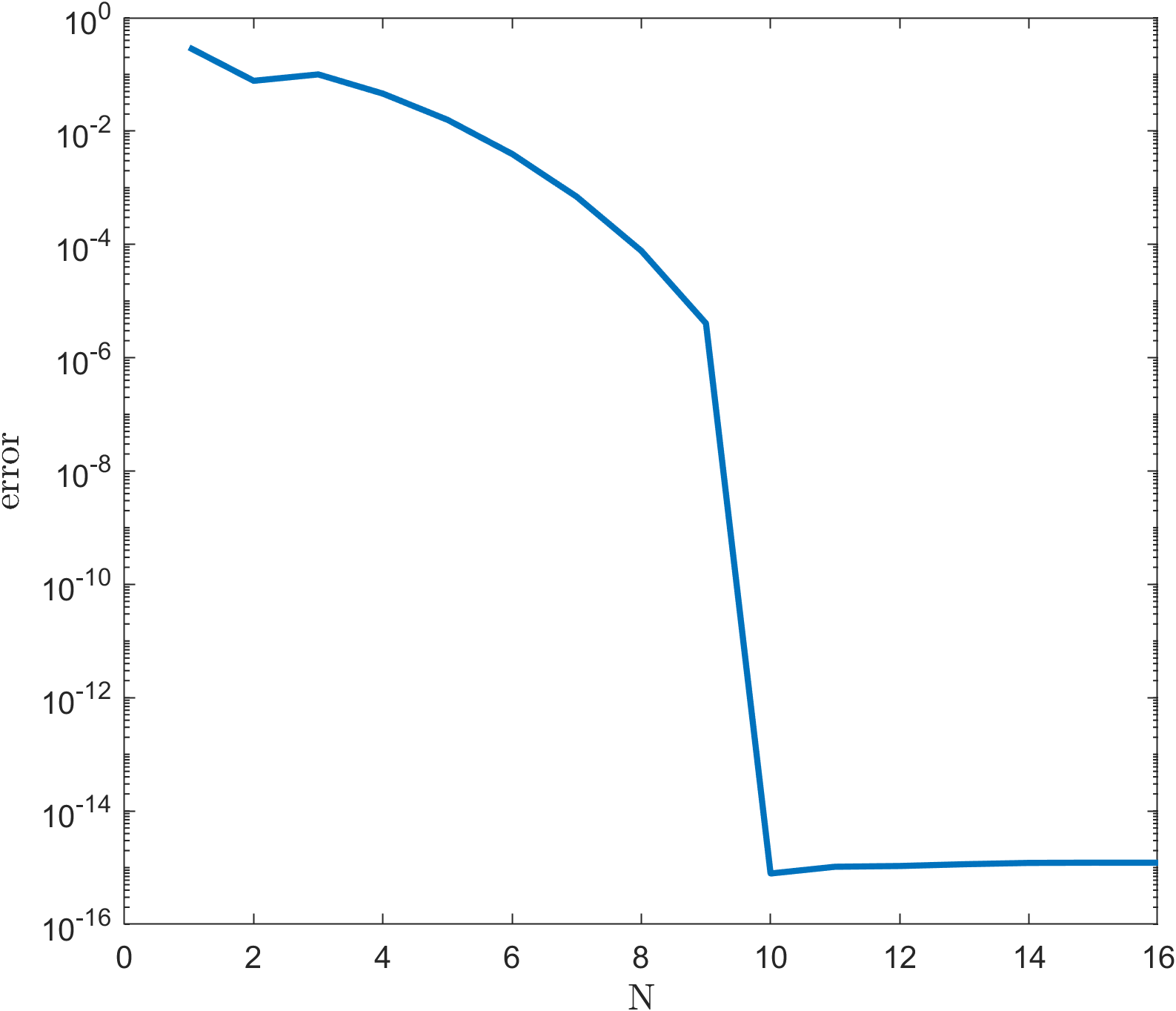}
\caption{Error in the numerical solution to \Cref{var} obtained by the JFP method.}\label{fig:ex_var}
\end{figure}

We now turn to an FIE that contains multiple FIOs of varying fractional orders, along with variable coefficients that appear either as prefactors to the FIOs or the solution $u(x)$
\begin{align}
u(x) + \sqrt{1+x}\mI^{1/3}u(x) + \mI^{1/2}[(1+\rotatebox{45}{\scalebox{0.5}{$\square$}})^{1/3}u(\rotatebox{45}{\scalebox{0.5}{$\square$}})](x) = f(x). \label{var}
\end{align}
The right-hand side is chosen as
\begin{align*}
f(x) = (1+x)^{3/2} + \left(\frac{\Gamma(5/2)}{\Gamma(17/6)} + \frac{\Gamma(17/6)}{\Gamma(10/3)}\right)(1+x)^{7/3},
\end{align*}
so that the exact solution is $u(x) = (1+x)^{3/2}$. Noting that the greatest common divisor of the fractional orders $1/3$ and $1/2$ is $1/6$, and that $f(x)$ is a polynomial in $(1+x)^{1/6}$, we seek the solution in the space spanned by $\mathbf{Q}^{0,1/6}$. The error in the numerical solution obtained using the JFP spectral method is shown in \Cref{fig:ex_var}. As expected, the error drops to machine precision at $N = 10$, since representing $u(x) = (1+x)^{3/2}$ requires only the basis functions $Q^{0,1/6}_j(x)$ for $j = 0, 1, \ldots, 9$.

\subsection{FDE boundary value problems}
\subsubsection{Basset--Boussinesq--Oseen equation}
As mentioned in the introduction, the matrix approximation to the FIO is also instrumental in solving FDEs. Our third example is the linear Basset--Boussinesq--Oseen equation 
\begin{align}
v'(t) +  \mD^{1/2}_t v(t) + v(t) = 0, \text{~~s.t.~~} v(0) = 1, \label{bbo}
\end{align}
where the fractional derivative operator $\mD^{1/2}_t$ is defined in the Caputo sense. 
The full Basset--Boussinesq--Oseen equation \cite{zhu} is a nonlinear FDE describing the motion and hydrodynamic forces on a small particle in an unsteady flow at low Reynolds numbers. Equation \Cref{bbo} is obtained from the full BBO equation by assuming that the particle position is known. By applying Laplace transform \cite{jag}, one can derive a closed-form solution to \Cref{bbo}:
\begin{align*}
v(t) = t^{-1/2}\left(\frac{1}{\sqrt{\pi}} - 2\Re\left(\frac{3+i\sqrt{3}}{6}E_{\frac{1}{2},\frac{1}{2}}\left(\frac{-1+i\sqrt{3}}{2}t^{1/2}\right)\right)\right).
\end{align*}

\begin{figure}[t!]
\centering
\subfloat[]{\includegraphics[width=0.48\linewidth]{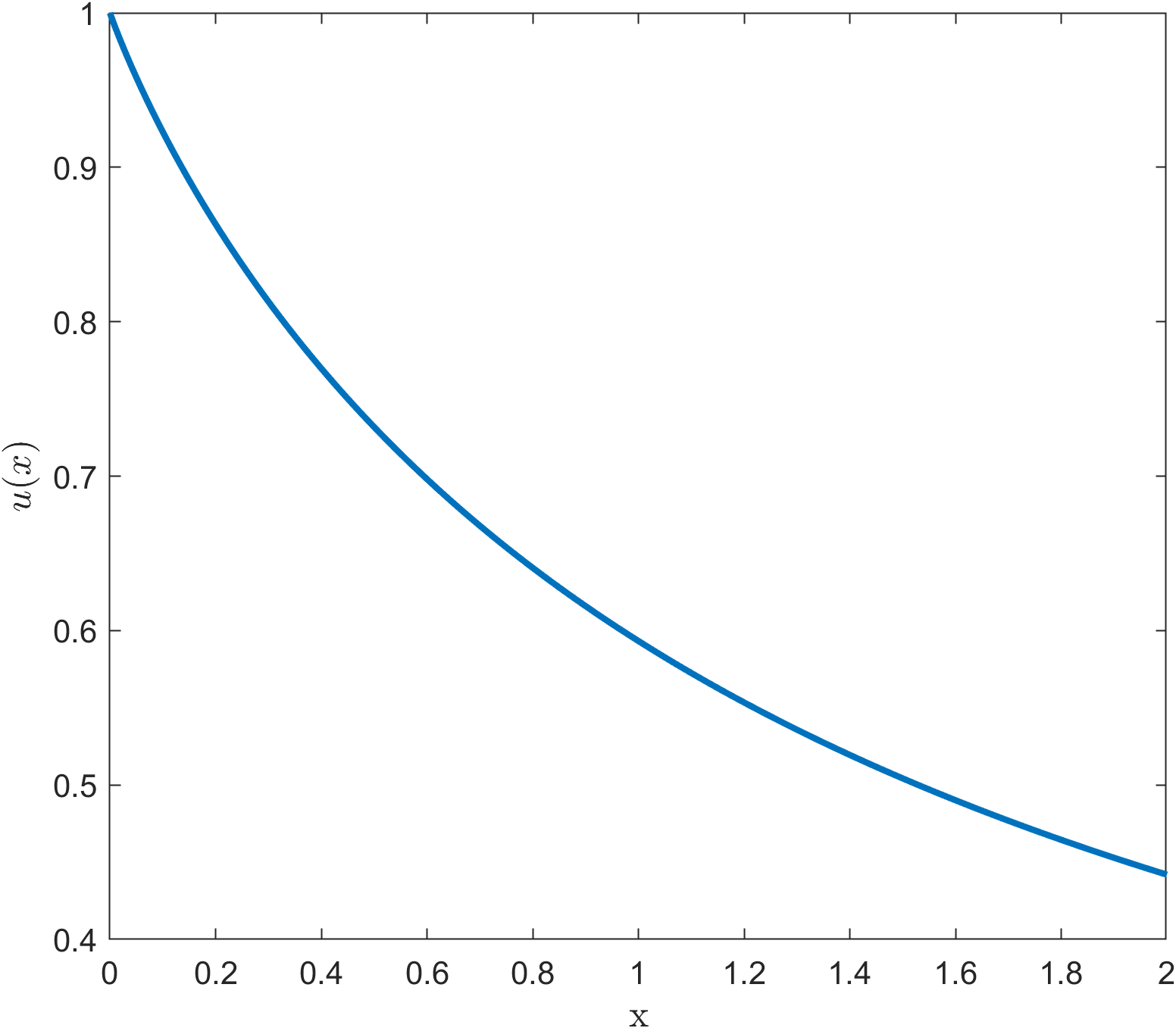}\label{fig:sol}}
\hfill
\subfloat[]{\includegraphics[width=0.48\linewidth]{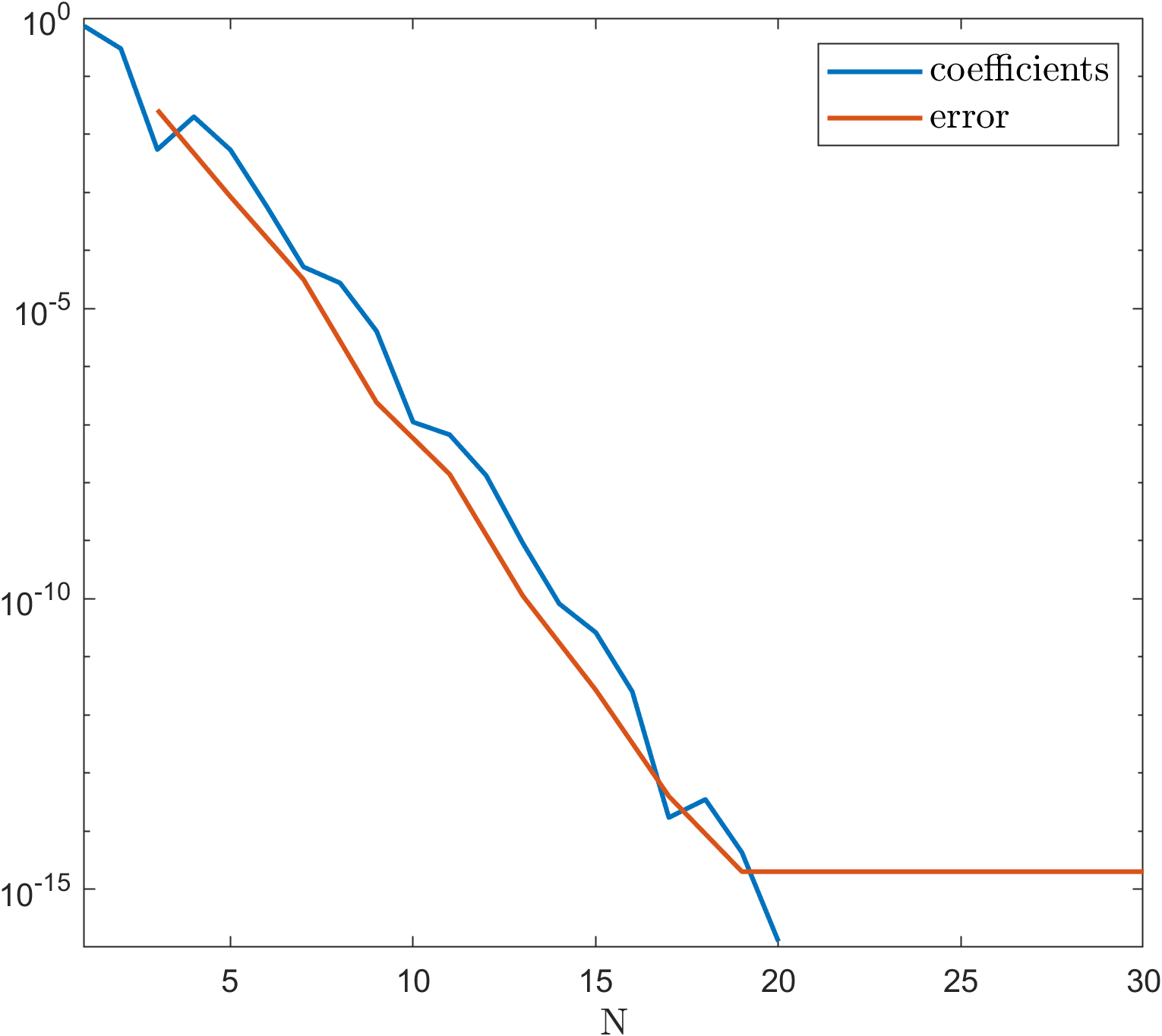}\label{fig:coeffs&error}}
\caption{(A) Numerical solution to \Cref{bbo}. (B) The JFP coefficients and the error of the computed solution.}
\label{fig:bas}
\end{figure}

Following the standard approach of integration reformulation, we let
\begin{align*}
v(t) = \mI^1 u(t) + 1,
\end{align*}
which transforms \Cref{bbo} into an FIE:
\begin{align*}
u(t) + \mI^{1/2} u(t) + \mI^1 u(t) + 1 = 0.
\end{align*}
The solution is sought as a series expansion in $\{Q_n^{0,1/2}(x)\}_{n=0}^{\infty}$. \Cref{fig:sol} shows the computed solution, while \Cref{fig:coeffs&error} illustrates the spectral decay of the JFP coefficients along with the error in the computed solution.

\subsubsection{Fractional Airy equation}
\begin{figure}[t!]
\centering
\subfloat[]{\includegraphics[width=0.3\linewidth]{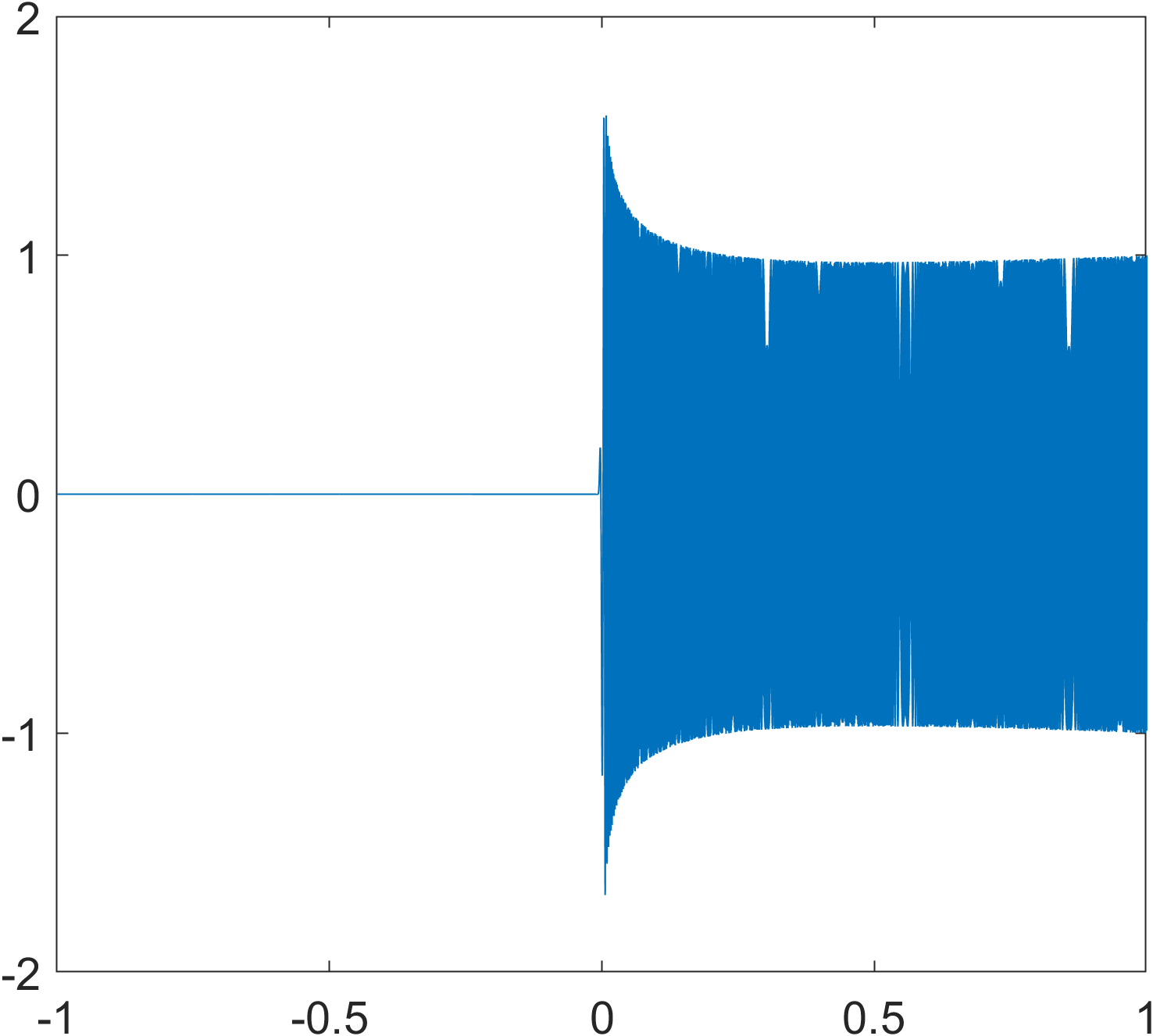}\label{fig:airy_real}}
\hfill
\subfloat[]{\includegraphics[width=0.3\linewidth]{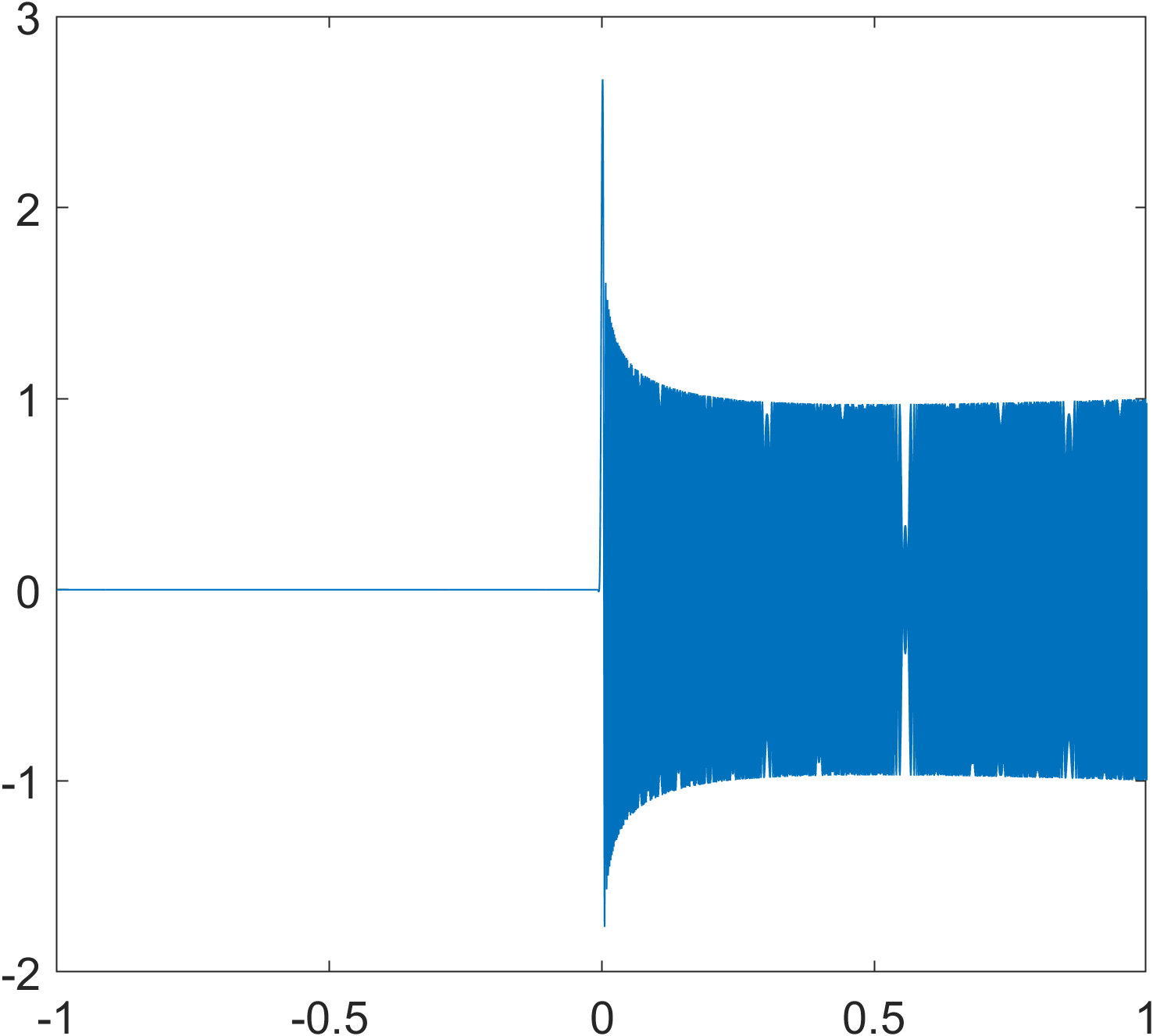}\label{fig:airy_imag}}
\hfill
\subfloat[]{\includegraphics[width=0.32\linewidth]{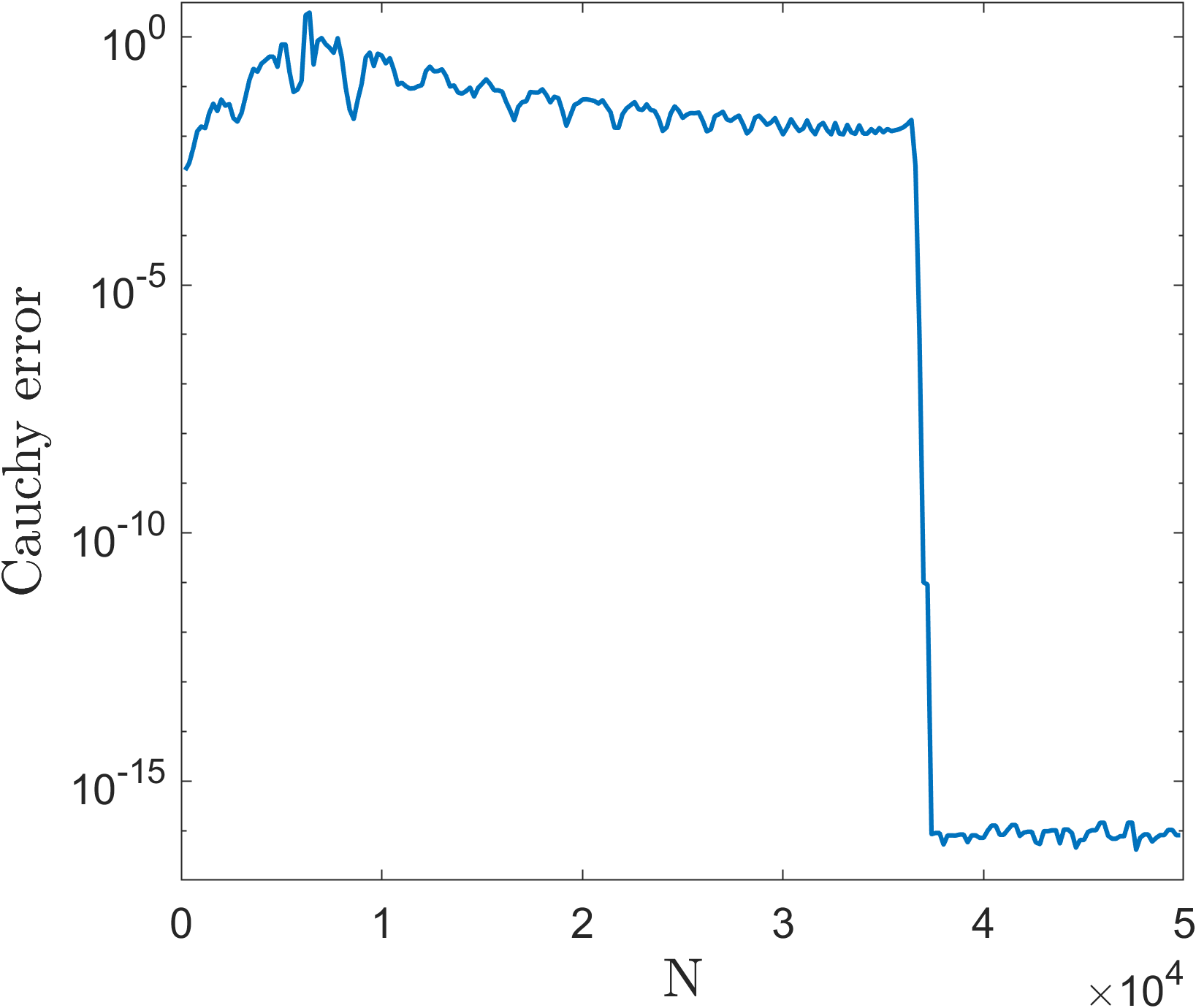}\label{fig:airy_cauchy}}
\caption{The (A) real and (B) imaginary parts of the numerical solution to \Cref{airy}. (C) Cauchy error.}
\label{fig:airy}
\end{figure}
Our second example for FDEs is the fractional Airy equation \cite{hal}
\begin{align}
\epsilon i^{3/2} \mD^{3/2}u(x) - xu(x) = 0, \quad x \in [-1,1], \quad \text{s.t. } u(-1) = 0,~~~ u(1) = 1. \label{airy}
\end{align}
We consider integral reformulation by 
\begin{align}
u(x) = \mI^{3/2}v(x) + a(1+x), \label{uv}
\end{align}
where $a$ is a constant to be determined. The ansatz in \Cref{uv} ensures that the boundary condition at $x = -1$ is satisfied automatically. Substituting \Cref{uv} into \Cref{airy} yields
\begin{subequations}
\begin{align}
\epsilon i^{3/2}\left(v(x) + \frac{a}{\Gamma(1/2)}\frac{1}{\sqrt{1+x}}\right) - x \mI^{3/2}v(x) + ax(1+x) = 0, \label{airyireqn}
\end{align}
subject to the boundary condition at the right endpoint:
\begin{align}
\mI^{3/2}v(1) + 2a = 1.
\end{align}\label{airyir}%
\end{subequations}

It follows from \Cref{airyireqn} that the solution and variable coefficients have expansions that contain a term of the form $1/\sqrt{1+x}$ and powers of $\sqrt{1+x}$. Hence, we represent them 
in the basis $\{ Q_n^{-1/2, 1/2}(x) \}_{n=0}^{\infty}$.
Let $\tilde{\mS}$ be the spectral approximation to the fractional integral operator $\mI^{3/2}$, and let $\hat{v}$ denote the coefficient vector of $v(x)$ in this basis. Replacing $\mI^{3/2}$ and $v(x)$ in \Cref{airyir} with $\tilde{\mS}$ and $\hat{v}$ gives the following linear system:
\begin{align*}
\left(
\begin{array}{c|c}
\mB \tilde{\mS} & 2 \\[1mm]
\hline\\[-2mm]
\epsilon i^{3/2} - \mM \tilde{\mS} & \hat{g} \\
\end{array}
\right)
\begin{pmatrix}
\hat{v}\\[2mm]
a
\end{pmatrix}
=
\begin{pmatrix}
1\\[2mm]
0
\end{pmatrix},
\end{align*}
where $\mB$ encodes the boundary condition (see \Cref{bc}), $\mM$ is the multiplication matrix (see \Cref{M}) with $c_j = 0$ for $j > 3$, and $\hat{g}$ is the $\{ Q_n^{-1/2, 1/2}(x) \}_{n=0}^{\infty}$ coefficient vector of 
\begin{align*}
g(x) = \epsilon i^{3/2}  \frac{1}{\Gamma(1/2)\sqrt{1+x}} - x(1+x).
\end{align*}

\Cref{fig:airy_real,fig:airy_imag} display the real and imaginary parts of the solution to \Cref{airy} for $\epsilon = 10^{-7}$, computed with a truncation size of $N = 50{,}000$. The corresponding Cauchy error, evaluated at each $N$ that is a multiple of $200$, using the numerical solutions at two consecutive values of $N$, is shown in \Cref{fig:airy_cauchy}. It can be observed that adequate resolution is achieved at approximately $N = 37{,}000$, beyond which the Cauchy error stabilizes around $10^{-16}$.

\subsection{FDE initial value problem}
The spectral approximations to FIOs are essential for achieving high-accuracy solutions to FDE initial value problems and time-dependent fractional PDEs via the spectral deferred correction (SDC) method \cite{dut, che3}. Consider the initial value problem
\begin{align*}
\mD^{\mu}_t u(t) = F(t, u(t)) ~\text{ s.t. }~ u(-1) = a,
\end{align*}
where $t \in [-1, 1]$ and $\mD_t^\mu$ denotes the Caputo fractional derivative. To apply the SDC method, we initialize with the trivial solution $u^0(t) = 0$ and evaluate the residual
\begin{align*}
\varepsilon (t) = a + \frac{1}{\Gamma(\mu)}\int_{-1}^{t}\frac{F\left(s, u^0(s)\right)}{(t-s)^{1-\mu}}\mathrm{d}s - u^0(t),
\end{align*}
where the FIO is approximated using the proposed spectral approximation. The correction term $\delta(t)$ is then computed by solving
\begin{align*}
\delta(t) = \frac{1}{\Gamma(\mu)} \int_{-1}^{t}\frac{F\left(s, u^0(t)+\delta(t)\right) - F\left(s,u^0(t)\right)}{(t-s)^{1-\mu}}\mathrm{d}s + \varepsilon(t),
\end{align*}
in which the fractional integral is approximated cheaply using a low-order finite difference scheme \cite{xu1} on a mapped Chebyshev grid
\begin{align*}
t_j = 2^{1 - 1/\beta}(1 + x_j)^{1/\beta} - 1, \quad x_j = \cos\left(\frac{j\pi}{N-1}\right),~~~ j = 0,1,\ldots,N-1,
\end{align*}
where $N$ denotes the total number of grid points. With the correction $\delta(t)$, the solution is updated, and the procedure is repeated until the residual $\varepsilon(t)$ falls below a prescribed tolerance.

\begin{figure}[t!]
\centering
\includegraphics[width=0.48\linewidth]{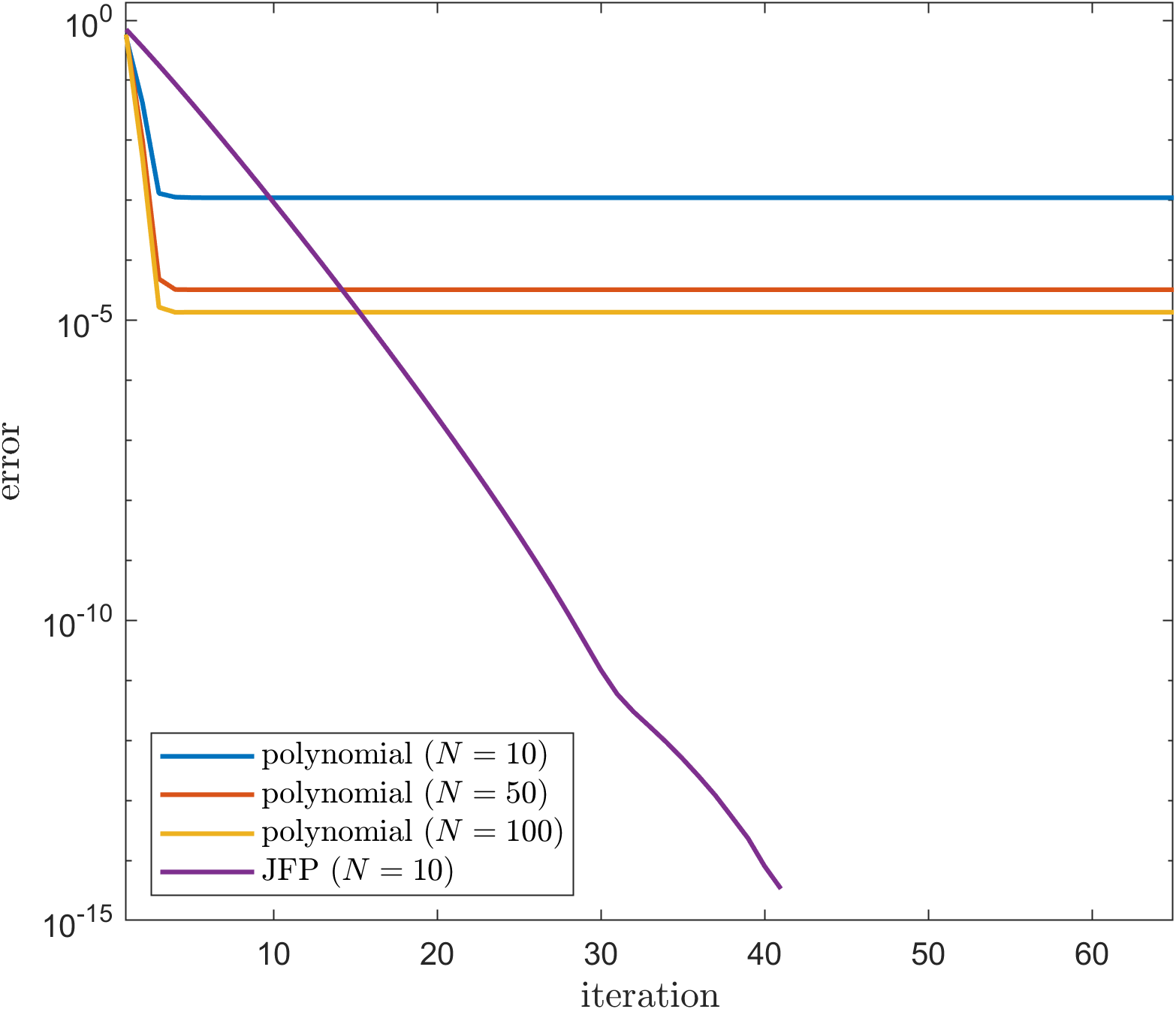}
\caption{Convergence of fractional SDC method based on polynomial and JFP bases.}\label{fig:ex_sdc}
\end{figure}

For illustration, we consider the case $\mu = 1/2$, $a = 0$, and
\begin{align}
F(t, u(t)) = u(t) + \sqrt{1 + t} - \frac{\Gamma(1/2)}{2}(1 + t), \label{F}
\end{align}
for which the exact solution is $u(t) = \Gamma(1/2)(1 + t)/2$. Accordingly, we employ the basis $\{Q_n^{0, 1/2}(x)\}_{n=0}^{\infty}$ with grid size $N = 10$. As demonstrated in \Cref{fig:ex_sdc}, the proposed implementation of the spectral deferred correction (SDC) method converges to machine precision within approximately 60 iterations. In contrast, the state-of-the-art SDC method described in \cite{che3} stagnates after only a few iterations. This breakdown in convergence is primarily due to the use of polynomial Lagrange interpolants and Legendre points in evaluating the residual $\varepsilon(t)$. For singular right-hand sides such as \Cref{F}, polynomial-based interpolants fail to provide sufficient accuracy, and the resulting inaccurate residual contaminates the correction steps, preventing further error reduction. Even when a large grid is employed, e.g., $N = 100$, a choice rarely adopted in practice, the accuracy plateaus at approximately $10^{-5}$.

\subsection{Fractional eigenvalue problem}
The matrix approximation to FIOs is also applicable to fractional eigenvalue problems. Suppose $\ell \geq 2$ is an integer, and that $\mu_1$ and $\mu_2$ satisfy $\ell - 1 < \mu_1 < \ell$ and $0 \leq \mu_2 < \ell - 1$, respectively, such that the ratio $\mu_2 / \mu_1$ is rational. Consider the eigenvalue problem
\begin{align}
    \begin{array}{l}
    -\mD_x^{\mu_1} u(x) = \lambda u(x) ~\text{ for }~ x \in [-1,1]\\[2mm]
     \text{s.t. }~ u^{(j)}(-1) = 0 \text{ for } j = 0,1,\cdots,\ell-2 ~\text{ and }~ \mD_x^{\mu_2}u(1) = 0,
    \end{array}\label{eigdif}
\end{align}
which plays an instrumental role in studying the properties of the two-parameter Mittag--Leffler function \cite{elo}. Here, $\mD_x^{\mu}$ denotes the Riemann--Liouville fractional derivative, and $u^{(j)} = \md^j u(x)/\md x^j$. Similar to FDEs, this problem can be reformulated as an FIO eigenvalue problem via integration:
\begin{align}
c (1+x)^{\mu_1 - 1} \mI^{\mu_1 - \mu_2}u(1) - \mI^{\mu_1}u(x) = \frac{1}{\lambda}u(x) , \label{eigint}
\end{align}
where the constant $c = \displaystyle \frac{\Gamma(\mu_1-\mu_2)}{\Gamma(\mu_1)}\mD_x^{\mu_2}[(1+\rotatebox{45}{\scalebox{0.5}{$\square$}})^{\mu_1 -1}](1)$. Let $\hat{v}$ and $\hat{u}$ denote the JFP coefficient vectors (infinite column vectors) of $(1+x)^{\mu_1 - 1}$ and the eigenfunction $u(x)$, respectively. Then the operator eigenvalue problem \Cref{eigint} can be expressed as
\begin{align}
\left( c \hat{v} \mB \mS^{\dagger} - \mS^{\ddagger} \right) \hat{u} = \frac{1}{\lambda} \hat{u}, \label{eigintmat}
\end{align}
where $\mB$ is again the infinite row vector of Dirichlet boundary condition \Cref{bc}, and $\mS^{\dagger}$ and $\mS^{\ddagger}$ are the matrix approximations to $\mI^{\mu_1 - \mu_2}$ and $\mI^{\mu_1}$, respectively.

\begin{figure}[t!]
\centering
\subfloat[]{\includegraphics[width=0.32\linewidth]{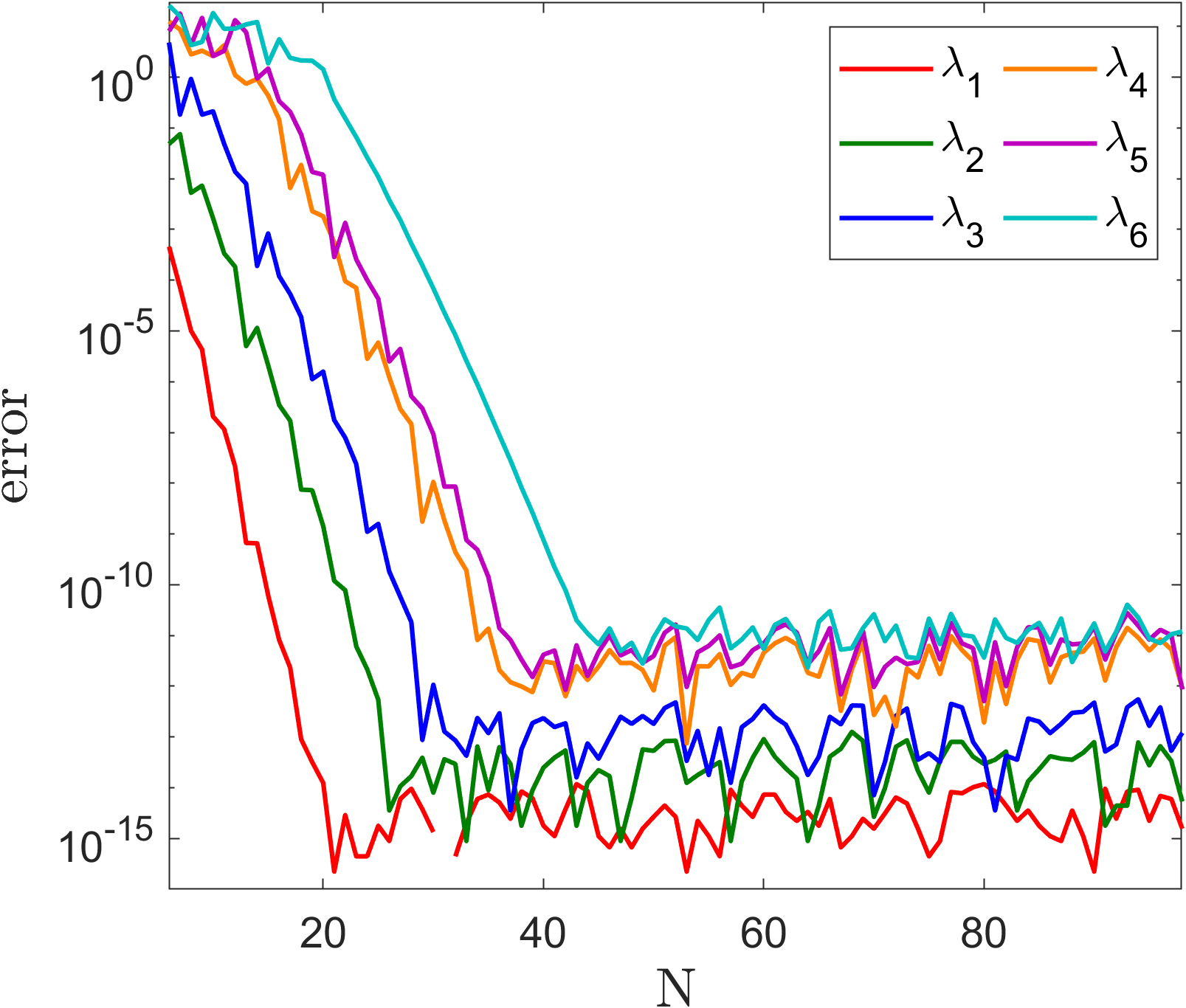}\label{fig:eig_cauchy}}
\hfill
\subfloat[]{\includegraphics[width=0.314\linewidth]{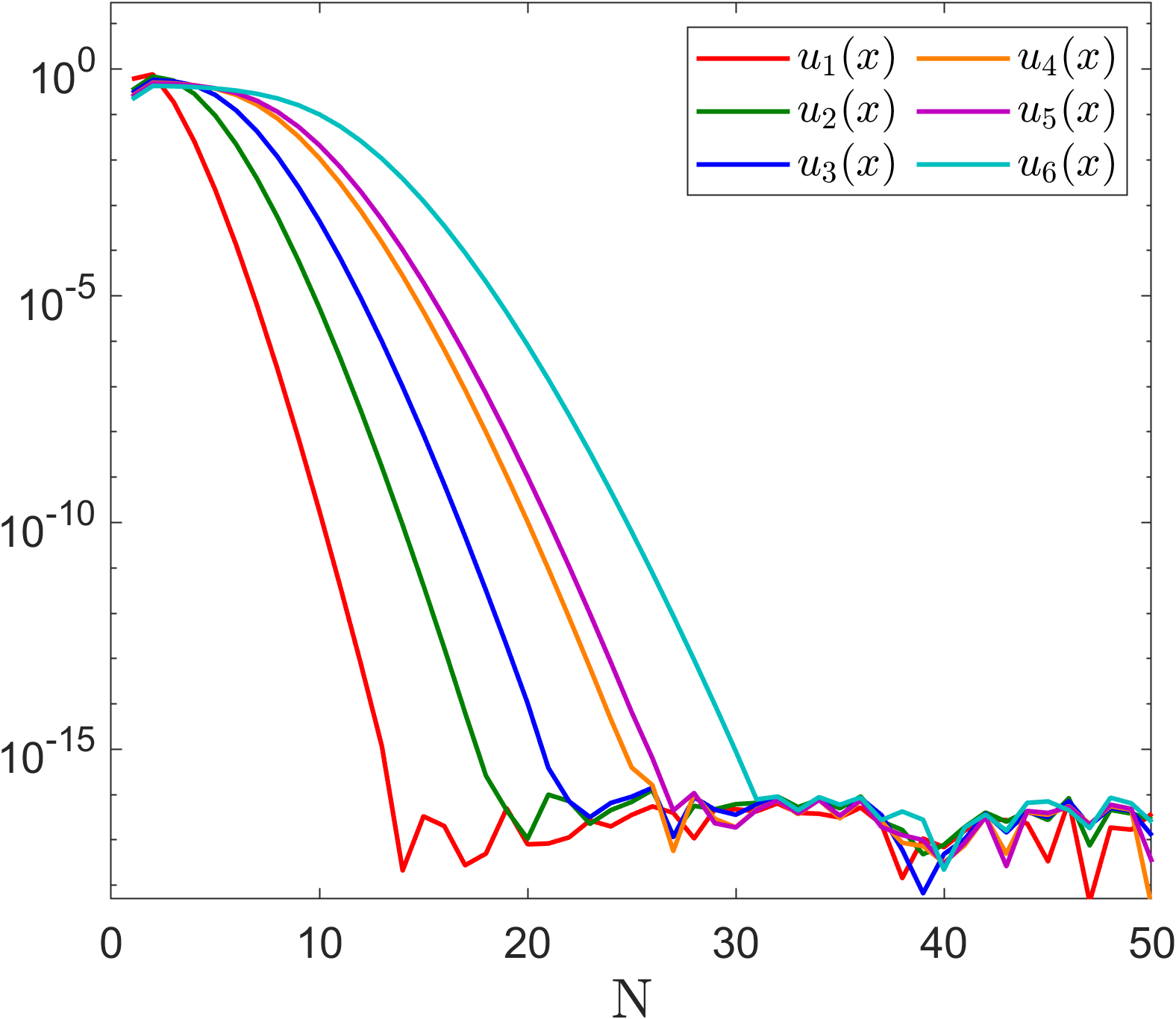}\label{fig:eigfun_coef}}
\hfill
\subfloat[]{\includegraphics[width=0.322\linewidth]{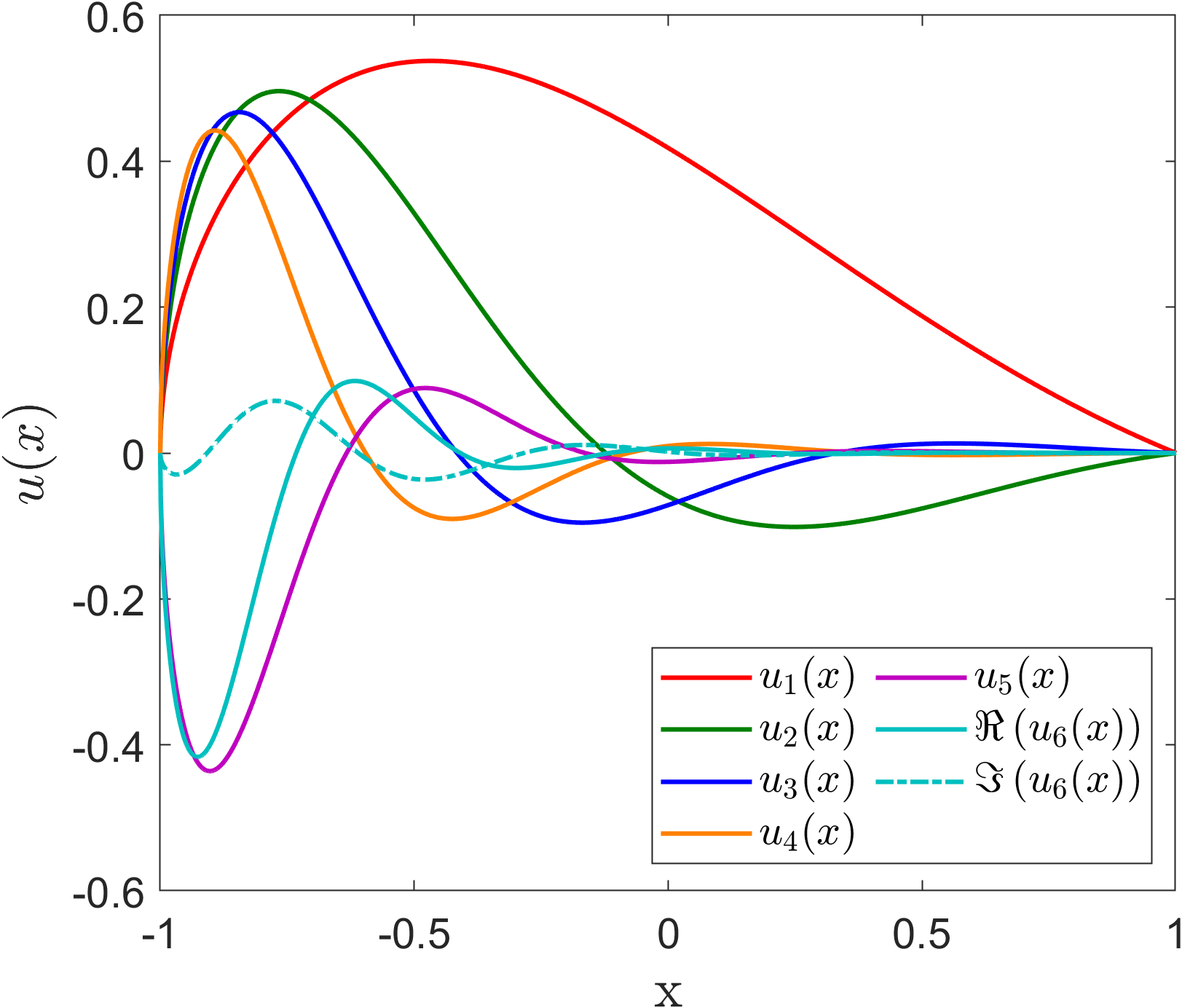}\label{fig:eigfun}}\\
\subfloat[]{\includegraphics[width=0.9\linewidth]{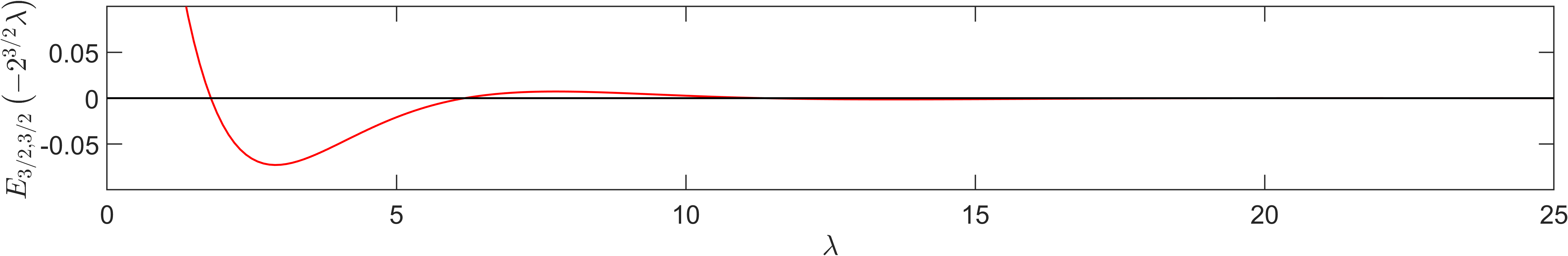}\label{fig:mlf}}
\caption{(A) Cauchy errors of the six eigenvalues of smallest modulus, obtained by solving the eigenproblem of consecutive truncation sizes using \texttt{eigs}. (B) JFP coefficients of the eigenfunctions. (C) Eigenfunctions corresponding to these six eigenvalues. (D) Plot of $E_{3/2,3/2}(-2^{3/2} \lambda)$ over $\lambda \in [0, 25]$.}\label{fig:eig}
\end{figure}

We set $\mu_1 = 3/2$ and $\mu_2 = 0$ in our experiment, and accordingly work with the basis $\{Q^{1/2,3/2}_n(x)\}_{n=0}^{\infty}$. To approximate the first six eigenpairs of the infinite-dimensional matrix eigenvalue problem \Cref{eigintmat}, we truncate the system and compute the six eigenvalues of smallest modulus, along with their corresponding eigenfunctions, using \textsc{Julia}'s \texttt{eigen}. This procedure is repeated with progressively larger truncation sizes. The Cauchy error of the computed eigenvalues is measured by the $2$-norm of the difference between eigenvalues obtained at two consecutive truncation sizes. Both the Cauchy errors and the computed eigenvectors, which approximate the JFP coefficients of the eigenfunctions, are examined for the formation of plateaus \cite{aur}. The iterative process is terminated only when both quantities exhibit plateaus (see \Cref{fig:eig_cauchy,fig:eigfun_coef}), which serves as an indicator of convergence and ensures accurate approximations to the eigenvalues and eigenfunctions \cite[Chap.~IV, \S 3.5]{kat}. The corresponding eigenfunctions are shown in \Cref{fig:eigfun}, with the real and imaginary parts of the sixth eigenfunction plotted separately.

We list the computed values of the six eigenvalues in \Cref{tab:valE}, and note that the first five are real, while the sixth is complex. Since $\lambda$ is an eigenvalue of \Cref{eigdif} if and only if it is a zero of $E_{\mu_1, \mu_1-\mu_2}(-2^{\mu_1}\lambda)$, we also include the norm of $E_{3/2,3/2}(-2^{3/2}\lambda)$ to indicate the accuracy of the computed values. The accuracy of the computed eigenvalues gradually deteriorates. This is not due to limitations of the JFP spectral method, but rather because the problem becomes increasingly ill-conditioned. This is illustrated in \Cref{fig:mlf}, where $f(\lambda) = E_{3/2,3/2}(-2^{3/2} \lambda)$ is plotted for $\lambda \in [0, 25]$. As shown, the function crosses the real axis at increasingly shallow angles, making the location of its zeros more sensitive to numerical errors.

\begin{table}[t!]
  \centering
  \renewcommand{\arraystretch}{1.4} 
  \setlength{\tabcolsep}{7pt}
  \caption{The six eigenvalues of smallest modulus, along with the corresponding values of $E_{3/2,3/2}(-2^{3/2}\lambda)$.}\label{tab:valE}
\begin{tabular}{c c c}
  \toprule
  index & eigenvalue $\lambda$ & $E_{3/2,3/2}(-2^{3/2}\lambda)$ \\
  \midrule
  1 & $1.794435495663993$ & $9.19\times 10^{-16}$ \\
  \midrule
  2 & $6.177290302782617$ & $6.31\times 10^{-16}$ \\
  \midrule
  3 & $11.359485354309392$ & $1.64\times 10^{-16}$ \\
  \midrule
  4 & $19.740438605284737$ & $-5.29\times 10^{-16}$ \\
  \midrule
  5 & $22.834767521795890$ & $-1.51\times 10^{-15}$ \\
  \midrule
  6 & $35.255579686924854\pm7.532188956823454i$ & $1.92\times 10^{-14}$ \\
  \bottomrule
\end{tabular}
\end{table}

\section{Conclusion and outlook}
The new algorithm that we propose is fast and stable in constructing the spectral approximation to FIOs. Such matrix approximations make the JFP spectral method practical for solving FIEs and FDEs and allow the fractional eigenproblems to be investigated numerically. The \textsc{Julia} implementation of this paper can be found at \cite{code}.

The method we propose shares some similarities with the approaches to constructing spectral approximations to the convolution operators of Volterra and Fredholm types \cite{xu,liu}. Preliminary results show that the recurrence-based method introduced in this paper, with some adaptation and specialized techniques, can be extended to virtually all linear operators. A unifying framework for operator approximation is an uncharted but exciting territory to be mapped in approximation theory.

\section*{Acknowledgments}
Many friends and colleagues have helped us along the way. In particular, we are grateful to Sheng Chen (BNU) for generously sharing his GLOFPG code, and to Xuejuan Chen (Jimei) for providing her implementation of the fractional SDC method. Tianyi Pu (Imperial) and Fanhai Zeng (Shandong) offered insightful feedback on an earlier version of the manuscript. We also thank our team members, Xi Chen and Kaining Deng, for their careful proofreading of the manuscript, as well as the two anonymous referees for their constructive comments.

\bibliographystyle{amsplain}
\bibliography{references}

\end{document}